\documentclass[a4paper]{amsart}

\usepackage[latin1]{inputenc}
\usepackage{amssymb,amsthm}
\usepackage[arrow,matrix,curve]{xy}

\usepackage{polyzetas, joe2004,joethm_en}
\usepackage{elliptique,connexions,perversDmod, 
            groupes_classiques, espace_modules}

\def\joenote#1{}

\let\s=\sigma
 
\def\ad{\text{\rm ad}}
\def\a{\mathbf t} 
\def\t{\mathbf t}

\def\alg{\text{alg}}
\def\aff{\text{aff}}
\def\talg{\t_\alg}
\def\salg{{\mathbf s}_\alg}

\def\odb{[\![} 
\def\cdb{]\!]} 	
\def\db#1{[\![\!{#1}]\!]} 

\title{Towards multiple elliptic polylogarithms}
\author{Andrey Levin}
\author{Georges Racinet}
\date{\today}

\begin{document}
\maketitle
\tableofcontents
\section{Introduction}
\ssect{History and motivation}
  \sss{Classical polylogarithms and $\Puntrois$} 
For the purposes of their proof of the weak version of Zagier's conjecture,
Beilinson and Deligne  \cite{BeilDel} gave a geometric interpretation
of the classical polylogarithm functions from the {\sc xviii}th century 
\begin{equation}\labeq{lisimple0}
        \Li_k(z)\ \ass\ \sum_{n>0} \frac{z^n}{n^k},\quad k\geq 1
\end{equation}
as periods of certain unipotent variations
of mixed Hodge structures on $\Puntrois$. 
On the differential side, the main feature is that the generating series 
$$ L(z)\ \ass\ \sum_{n\geq 0} \Li_n(z)t^n$$
satisfies the following differential equation
$$ dL(z)\ =\ E_0\dlog(z) + E_1\dlog(1-z), $$
where $E_0, E_1$ are explicit topologically nilpotent linear
operators on $\serc{\CM}{t}$. 

\sss{Generalisations}
One may consider the following generalisation with several integer indices:
\begin{equation}\labeq{Li0}
\Li_{s_1,\ldots,
  s_r}(z)\ \ass\ \sum_{n_1>n_2>\cdots>n_r>0}\frac{z^{n_1}}{n_1^{s_1}\cdots
  n_r^{s_r}}, \quad r>0, s_i>0  
\end{equation}
The noncommutative generating series 
  $$1+\sum_{s_1,\ldots,s_r
  >0} \Li_{s_1,\ldots,s_r}(z)e_0^{s_1-1}e_1\cdots e_0^{s_r-1}e_1$$ 
  is a horizontal section of the trivial vector bundle $\PC$ on
  $\Puntrois$ with fiber $\sernc{\CM}{e_0,e_1}$, the completed free
  associative algebra with generators $e_0,e_1$, and connection 
  \begin{equation}\labeq{KZ3red}
     \nabla = d - e_0\dlog(z)-e_1\dlog(1-z),
  \end{equation}
   where $e_0, e_1$ act by left
  multiplication within $\sernc{\CM}{e_0,e_1}$. 

  As a vector bundle with
  connection, $\PC$ is universal in the sense that any nilpotent
  vector bundle with regular singular connection on $\Puntrois$ 
  is isomorphic to 
  the trivial bundle with fiber $V$ and connection given
  by \refeq{KZ3red}, for some $E_0,E_1\in\GL_\CM(V)$. This
  is best expressed in the language of 
  tannakian categories.  In the terminology of Deligne \cite{DelPi1}, 
  $\PC$ is the De Rham version of
  a \emph{fundamental torsor} of motivic paths on $\Puntrois$
  (loc. cit., \S 12). This
  fundamental torsor provides the geometric interpretation of these
  generalized polylogarithms. 

  One can push the definition further and consider the \emph{multiple
  polylogarithm} 
\begin{equation}\labeq{multpol}
\Li_{s_1,\ldots,
  s_r}(z_1,\ldots,z_r)\ \ass\ \sum_{n_1>n_2>\cdots>n_r>0}\frac{z^{n_1}\cdots
  z^{n_r}}{n_1^{s_1}\cdots
  n_r^{s_r}}  
\end{equation}
  which can be interpretated in terms of fundamental groups and torsors of 
  the moduli spaces $\MZ{n}$ of genus $0$ curves with $n$ marked points
  \cite{Gonch2001}. 

\sss{}
  Special values of these functions, and especially the \emph{multiple
  zeta values} (MZVs):
\begin{equation}\labeq{MZV}
        \z(s_1.\ldots,s_r)\ \ass\ \sum_{n_1>\cdots>n_r>0}\frac{1}{n_1^{s_1}\cdots n_r^{s_r}}
\end{equation}
 are of prime importance in the theory of mixed
  Tate motives. There is a very precise conjectural
  picture \cite{DelGon, Gonch98} of the algebraic rational relations
  that occur among them. In particular, the pro-unipotent
  Grothendieck-Teichmüller 
  program in genus $0$  relates them to a conjectural
  description of nilpotent quotients of 
  $\text{\rm Gal}(\ov\QM/\QM(\mu_\infty))$ \cite{IhICM}. 
  It may be worthwile to notice that one can extract from
  explicit formulas like \refeq{MZV} a lot of information in a combinatorial
  way \cite{joeIHES, ENR}, whose interpretation in terms of
  differential equations is rather non-trivial \cite{Gonch2001}. 

  \sss{Elliptic polylog}Another direction of generalization was taken by
  Be\u\i linson and the first author \cite{BeilLev}. The polylogarithmic sheaf 
  can be characterized
  canonically by its formal properties with respect to the inclusion of
  $\Puntrois$ into $\Gm$ \cite{BeilDel}. This abstract description 
  has an elliptic parallel that allows to define elliptic polygarithmic 
  mixed sheaf on a punctured elliptic curve. There is a relative version for
  families of elliptic curves, and a continuity property to the ordinary 
  polylogarithmic sheaf at cusps.

  This approach was further investigated by Wildeshaus \cite{Wildes} who 
  generalized it to the complement of a mixed Shimura variety 
  into another.
  \sss{Our goal: multiple elliptic polylogarithms}In this article, we start 
  to study multiple elliptic polylogarithms in terms of De Rham 
  fundamental torsors approach of \cite{DelPi1}. 

  It seems rather natural that the sheaf theoretic elliptic
  analog of \refeq{Li0} should come from the whole fundamental torsor of
  paths on a 
  punctured elliptic curve, whereas the full multiple elliptic
  polylogarithms should be related to fundamental torsors of the moduli
  stacks $\M{1,n}$.

  At this point, one should mention that in \cite{DelPi1}, Deligne used 
  systems of realisations (among them the De Rham realisation) 
  as a substitute for the conjectural category of mixed Tate motives. Such a
  category is now well defined, and the motivic fundamental groups and torsor 
  have received direct motivic definitions \cite{DelGon}. 
  However, as of this writing, mixed elliptic motives seem to be still far
  out of reach. 

\ssect{What we do}
  We provide here an explicit description of the De Rham fundamental
  torsor of an elliptic curve $X$ minus its origin and extend this to
  families of those, \ie, to $\M{1,2}$.
  In equivalent and more down-to-earth words, we give a
  complete and explicit classification of vector bundles with
  (relatively) nilpotent connection on a (family of) elliptic curves.  

  Explicit formulas in the genus $0$ case, such as \refeq{KZ3red}, 
  arise from the fact that any nilpotent connection on $\Puntrois$
  is the restriction of a connection with simple poles at $0,1,\infty$
  on a \emph{trivial} bundle on $\PM^1$. In general, the relevant category to 
  consider for the De Rham fundamental group(oid) of a variety given as 
  complement of a normal crossing divisor $D$ in a proper variety $X$ is
  the category $\nilpconn(X; D)$ of nilpotent vector bundles 
  \emph{on the whole of $X$} with connection having simple poles along $D$ 
  \cite{Del_singreg, BoAl, DelPi1}:
    \begin{equation}\labeq{good_case}
      \nabla\colon \VC\ \longto\ \VC\ot\Omega^1_X(D)
    \end{equation}
  Such nilpotent bundles are trivial if one assumes that $H^1(X, \OC_X) = 0$.
  We refer informally to this situation as 
  \emph{Deligne's good case} \cite[12.1]{DelPi1}. A fair amount of 
  \cite{DelPi1} relies on this hypothesis.

In this article, we circumvent the fact that an elliptic curve does
not satisty \refeq{good_case} by using
the analytic 
description of elliptic curves as quotients of type $X=\CM^*/q^\ZM$ and 
pulling everything back as trivial bundles on the affine $\CM^*$, with 
$q$-action. The correspondence between 
$q$-difference equations and vector bundles on elliptic curves is very
classical, and becomes quite simple in the nilpotent case, 
as well as the condition
that a connection on $\CM^*$ has to fulfill to be a pullback from
$X$. 

\ssect{Plan of the article} 

In section \ref{sec_fw}, we provide the description of the
fundamental De Rham torsor and fundamental group 
for a given punctured elliptic curve
$X=\CM/\Lambda\setminus e$, by means of a full analytic description of the
category of vector bundles with meromorphic connection on $X$ and 
simple poles at the unit $e$. This part is very elementary. 

The end result is as close as one could hope to 
Deligne's good case: we obtain a canonical abstract
basepoint (ie., a fiber functor) $\omega_\dr$,
at which the fundamental Hopf algebra is canonically identified with
the free Lie algebra in two 
generators. Furthermore, the universal torsor $\PC$ at $\omega_\dr$ and its 
connection are completely explicit. This could be thought of
as the elliptic analog of \refeq{KZ3red}.

In section \ref{sec_univ}, we define a vector bundle with flat
connection $\PG$ on the standard family $\EM$ of elliptic curves over
Poincaré's upper-half plane, such that the fiber at any $\tau$ is the former
$\PC$ of the corresponding elliptic curve. 
Furthermore, $\PG$ is $\SL_2(\ZM)$-equivariant. 

In section \ref{sec_rel}, we turn to arbitrary smooth families of
elliptic curves and to relatively nilpotent connections with
simple poles at the unit section. For this, we merge the
results of the previous sections: $\PG$ descends to a vector bundle
$\PG_{X/S}$ on each family $X\to S$. The main result, theorem
\ref{thm_rel}, then mostly says that $\PG_{X/S}$ 
is the relative fundamental torsor at
some abstract fixed fiber functor $\omega_\dr$, in the sense of 
\cite{spitzdiplom}.  
In short, this is an explicit analytic description of the De Rham torsor 
of paths on the moduli stack $\MC_{1,2}$, relative to $\MC_{1,1}$.

We then treat the case of geometric base points, \ie{}, sections of our
family of punctured elliptic curves and derive a differential
equation for parallel transports along families of paths have to fulfill. 
We believe this to be also a nice and concrete illustration of basic 
2-categorical aspects of the relative tannakian theory.

In the last section, we exhibit natural $\QM$-structures on the 
previously constructed analytic fundamental torsors and their connections. 
The construction goes mostly through the rewriting of the previously obtained 
formulas in terms of elliptic functions and modular forms.

\ssect{Conventions and notations}
Wherever possible without leading to confusion, we'll use the word
``connection'' to actually mean a vector bundle with connection. 

The letter $X$ will usually denote an elliptic curve, or a smooth family of
these on a base $S$. Since we use alternatively additive and multiplicative
conventions, we'll always denote the unit (section) by $e$. 

Straight letters usually refer to objects living over a point: vector spaces, 
algebras, a typical example would be the algebra $R$ from \ref{nota_R}.
We'll use curly letters for objects living over a space, for instance the 
fundamental torsor $\PC$ of \ref{fw_torsor} or the algebra bundle with 
connection $\RC$ of \ref{def_Rrel}. We'll use capital gothic letters for 
sheaves over a family, like the $\PG$ from \ref{def_PG}.

Exceptions to these conventions are Poincaré's upper half-plane $\HM$, the 
family on elliptic curves $\EM\to\HM$. 

\sss{Eisenstein functions and series}
We recall the notations from \cite[chap. III \& IV]{WeilEll} that we'll 
throughout this article.

Consider a lattice in $\CM$ generated by $\tau$ and $1$, with $\Im(\tau)>0$. 
We take $\xi$ as coordinate on $\Bbb C$ and set 
$z=\exp(2\pi i \xi)$, $q=\exp(2\pi i \tau)$

The symbol $\sum'$ denotes the summation over nonzero elements of a lattice
and $\sum'_e$ denotes some variant of regularization of the divergent sums,
known as {\em Eisenstein summation}. With this, we set:
$$ e_k(\tau)=\sum'_e(n\tau +m)^{-k} \qtext{and} E_k(\xi,\tau)=
\sum'_e(\xi+n\tau +m)^{-k}$$
The $e_k$ are modular functions for $k\ge 4$, and
$E_k$ are elliptic functions for $k\ge 2$. Moreover one can express the
Weierstrass elliptic functions and modular forms as
$\wp=E_2-e_2$, $\wp'=-2E_3$, $g_2=60e_4$ and $g_3=140e_6$.
It is known that the $e_k$ are rational polynomials in $e_4$ and $e_6$ for
$k\ge 8$, while the $E_k$ are rational polynomials in 
$E_2-e_2$, $E_3$, $e_4$ and $e_6$ 
for $k\ge 4$. Moerover the rings $\Bbb Q[e_k]$ and $\Bbb Q[E_k,e_k]$ are  
stable under the derivation $(2\pi i)\partial_\tau$.

The basic theta function is defined by the Jacobi product formula
$$\theta(\xi, \tau)=-iq^{1/8}(z^{1/2}-z{-1/2})\prod_{j=0}^\infty 
(1-zq^j)(1-z^{-1}q^j)(1-q^j).$$

We have $E_1=\partial_\xi\theta$.
\sss{Lie algebras} 
Recall the notation for the adjoint action of a Lie algebra on itself 
$\ad_a=(x \mapsto [a,x])$.


\section{The De Rham fundamental groupoid of a punctured elliptic curve}\label{sec_fw}

In this section, we consider an elliptic curve $X=\CM/\Lambda$, where
$\Lambda$ is the lattice of $\CM$ generated by $\tau$ and $1$, for
some fixed $\tau\in\HM$. As usual, we'll denote by $\xi$ the coordinate 
on $\CM$ and $z=\exp(2\pi i\xi)$ the coordinate on $\CM^*$.
Since we'll switch frequently between this
additive description and the multiplicative $X=\CM^*/q^\ZM$, we'll
denote the unit of $X$ by $e$. 
  
\ssect{Nilpotent vector bundles}\label{ss_bdl_fw}
Before starting to investigate connections on $X$, we have to describe
in analytic terms the nilpotent vector bundles themselves. Needless to
say, apart perhaps from the formulation, there is definitely nothing original 
in this subsection.

\sss{}
Let $V$ be a $\CM$-vector space and $m\in\endo{}{V}$. 
It is a well known fact that the equation
\begin{equation}\labeq{diff}
v(qz) = mv(z),
\end{equation}
%
where $m\in\GL(V)$, gives rise to a vector bundle $\VC$ on $X$, whose
sections over an open 
subvariety $U$ are the analytic solutions of \refeq{diff} defined over the 
preimage
of $U$ in $\CM^*$. We'll frequently call $\VC$ the \emph{bundle with
multiplier $m$}. An equivalent considerations can be made in the additive 
situation, \ie, for $X=\CM/\Lambda$ that we shall use liberally when it's 
more convenient, for instance for modular considerations.

\sss{}The pairs $(V,m)$ form a category in the usual way: an arrow from 
$(V_1,m_1)$ to $(V_2,m_2)$ is simply a linear map $f$ which intertwines 
$m_1$ and $m_2$, \ie, statisfies $fm_1=m_2f$. 
It can be further equipped with a tensor structure by the standard rule:
$$(V_1,m_1)\ot(V_2,m_2)\ \ass\ (V_1\ot V_2, m_1\ot m_2)$$
and inner homomorphisms:
$$\inthom{}{(V,m_1)}{(V,m_2)}\ \ass\ (\hom{\vects{\CM}}{V_1}{V_2},
f \mapsto m_2fm_1^{-1})$$ 
Altogether, we have a rigid tensor category, which is nothing but 
the category $\reps\ZM$ of representations of $\ZM$. 

The trivial object is simply $(\CM, \id)$. Therefore the
nilpotent objects are the pairs $(V,m)$ such that $m$ is unipotent.    

Given two objects $(V_1,m_1)$ and $(V_2,m_2)$, any linear map $u$ that
interwines $m_1$ and $m_2$ maps obviously solutions of \refeq{diff} to
solutions \refeq{diff}, hence our construction of bundles extends to 
a functor $\GC$.

\begin{prop}\label{thm_bdl}
	The functor $\GC\colon (V, m)\mapsto \VC$ is an
        equivalence of categories
  	 $$ (\text{\rm vector spaces with unipotent
	  endomorphism})\longto\nilpbdl(X)$$ 	

In particular, $\nilpbdl(X)$ is a sub-abelian category of the category
of coherent $\OC_X$-modules and is tannakian. 
\end{prop}

\sss{Global sections} Before proceeding to the proof of the above
theorem, we gather here some very elementary yet useful properties of global
sections of nilpotent bundles.   
\begin{prop}\label{prop_hol_sect}
	Let $m\in\GL(V)$ be unipotent. The global holomorphic sections of the
	bundle with multiplier $m$ are the constants fixed by $m$. 
\end{prop}
\begin{proof}
Thanks to the functoriallity of $\GC$, it is enough to prove this for
$(V,m)=(\CM^n,\id+J_n)$, where $J_n$ is the 
lower Jordan block of rank $n$. We prove by induction on $n$ that, for
any global section $s=(s_1,\ldots,s_n)$, we have $s_j=0$ for $j<n$ and
$s_n$ is a constant. 

The equation \refeq{diff} can be expressed as:
\begin{equation}\labeq{diff_jordan}
	s_{j+1}(qz) = s_j(z) + s_{j+1}(z),\qtext{for} j=0,\ldots,n-1,
\end{equation}
where we put $s_0=0$.

For $j=0$, this equation means that $s_1$ is an entire elliptic
function, hence must be constant. If $n=1$, we are done. Otherwise, for
$j=1$, the set of meromorphic 
solutions of \refeq{diff_jordan} is easily seen to be 
$\{ (2\pi i)^{-1}s_1E_1 + \phi \}$,
where $\phi$ runs over elliptic functions. Since $E_1$ has 
simple poles at lattice points only, which can't be the case for $\phi$, we
must have $s_1=0$ for $s_2$ to be entire. 
Now $(s_2,\ldots,s_n)$ is a global section of $\GC(\CM^{n-1}, J_{n-1})$.
\end{proof}
We'll have more solutions if we allow logarithmic poles, but there is
some rigidity which will prove useful in the course of the proof of theorem
\ref{thm_equiv_fw}: 
\begin{prop}\label{prop_unicite_sect}
     Let $m$ be a unipotent endomorphism of $V$. There is at most
     one meromorphic section of the bundle with multiplier $m$ having a
     simple pole at $e$ and given degree $0$ 
     component in its Laurent series expansion near $e$.
\end{prop}
\begin{proof}
   It's enough to prove that any such section $s$ with
   zero constant term vanishes. 
   Let $s_1,\ldots,s_n$ be the coordinates of $s$ in a basis such that the
   matrix of the multiplier is lower triangular.

   Again, we have $s_1(qz)=s_1(z)$, hence $s_1$ is elliptic. Since it
   has at most one simple pole, it must be constant, hence identically
   zero, by the assumption on the constant term. Now applying $s_1=0$ in our 
   triangular system of $q$-difference equations yields the induction.
\end{proof}

\begin{proof}[Proof of proposition \ref{thm_bdl}] 
	The unit object of the source category is $(\CM, \id)$. Its
	image is the sheaf of local elliptic functions, \ie, the trivial
	bundle $\OC_X$. On the other hand, $\GC$ is exact, since
	exactness can be checked on fibres, and $\GC_{|z}$ is simply the 
        forgetting of multipliers.

	\noindent{Tensor structure} 
	
	\noindent\emph{Fully faithfullness}
	The faithfullness is obvious, and can be checked on any fibre.
	For fullness, we have to prove that any morphism
	$\GC(V_1,m_1)\to\GC(V_2,m_2)$ comes from a constant $u:V_1\to
	V_2$ such that $u=m_2 u m_1^{-1}$. Since homomorphisms of
	bundles are global sections of $\homint{}{\VC_1,\VC_2}$, this
	follows from proposition \ref{prop_hol_sect}.

	\noindent\emph{Essential surjectivity} By Atiyah's theorem 5
	of \cite{Atiyah_fibres}, there is up to isomorphism exactly one
	indecomposable nilpotent bundle on $X$ of rank $n$. Using the fully
	faithfullness, we see that 
	the bundle associated to $\CM^n$ 
	and $m=\id+J_n$, where $J_n$ is the Jordan block of rank $n$,
	has those required properties. 
\end{proof}
\sss{Remark} Let $S$ be a space (in some sense) with action of a group $G$ and 
assume $S/G$ makes sense in such a way that  $p\colon S\to S/G$ is a 
$G$-bundle. he pullback by $p$ induces an equivalence from the category of 
bundles on $S/G$ to the category of $G$-equivariant bundles on $S$. 
In our case, where $S$ is the 
multiplicative group, any vector bundle is trivial. The multiplier 
can be thought as the expression of the incompatibility between an arbitrary 
trivialisation and the $G$ action. The additional information we obtained is 
that \emph{constant} multipliers are enough.

\ssect{Nilpotent connections with simple poles at the origin}
\sss{A two variable Jacobi form}
The following function was introduced by Kronecker \cite{KroneckerF},
rediscovered by Zagier \cite{ZagModTheta}, and considered by the
first author in the context of elliptic polylogarithms \cite{LevinAnal}:
\begin{eqnarray*}
F(\xi , \alpha; \tau )\ \ass \ (2\pi i)\left(
1 - \frac {1}{1-z} -\frac{1}{1-w} - \sum_{m,n=1}^{\infty}
(z^m w^n -z^{-m} w^{-n}) q^{mn}\right),
\end{eqnarray*}
where $q=\exp (2\pi i \tau ),\quad
z=\exp (2\pi i \xi ),\quad w=\exp (2\pi i \alpha )$.\joenote{Copypaste
from one paper of Andrey for notations to be put before this}
We recall here some of the statements of \cite[3. Theorem]{ZagModTheta}

The double series $F$ converges in the domain $\{\Im \tau >\Im \xi >0\,
\Im \tau > \Im \alpha > 0\}$ and extends meromorphically to all
values of $\xi$ and $\alpha$. It has simple poles at divisors $\xi
=m+n\tau$ and $\alpha = m' + n'\tau$ and can be expressed by means of
the Jacobi theta function:\footnote{This is actually the definition Kronecker
started with.} 
\begin{equation}
F(\xi ,\alpha ;\tau )=\frac{\theta'(0;\tau ) \theta(\xi + \alpha ;\tau )}
{  \theta (\xi ;\tau ) \theta (\alpha ;\tau )}. \label{theta}
\end{equation}
The residue of $F$ at $\xi=0$ is $\alpha^{-1}$. 

The elliptic and modular properties of $F$ are
rather nice:  
\begin{eqnarray}
F(\xi+1 ,\alpha ;\tau )&=&F(\xi ,\alpha ;\tau );\labeq{F_per}\\
F(\xi +\tau ,\alpha ;\tau )&=
&\exp (-2 \pi i \alpha )F(\xi ,\alpha ;\tau );\labeq{F_quasiper} \\
F(\frac{\xi}{c\tau +d} ,\frac{ \alpha}{c\tau +d} ,
\frac{a\tau +b}{c\tau +d} )&=&
(c\tau +d)\exp (2\pi i\frac{c\xi \alpha }
{c\tau +d})F(\xi , \alpha ;\tau ).\label{dermod}\labeq{F_mod}
\end{eqnarray}

$F$ can also be expressed as the exponential of the generating series
in one variable of Eisenstein functions of the other variable:
\begin{equation}\labeq{F_Eis}
   F(\xi,\alpha;\tau) = 
	\exp\left(-\sum_{k\geq 1} \frac{(-1)^k\xi^k}{k}(E_k(\alpha; \tau)-e_k(\tau)\right)
\end{equation}
This statement can easily be deduced from Zagier's ``logarithmic
formula'' for $F$ (loc. cit., {\it(viii)}) and the power series
expansion formula for $E_n$ \cite[III, (10)]{WeilEll}

\sss{}Although they won't be needed before section \ref{sec_univ}, we
gather here further properties of $F$, also easily deduced from
classical properties of theta series.  

It satisfies the ``mixed heat'' equation:
\begin{equation}\labeq{mixed_heat}
      2\pi i\frac{\partial F(\xi,\alpha;\tau)}{\partial\tau}
      = \frac{\partial^2F(\xi,\alpha;\tau)}{\partial\xi\partial\alpha}
\end{equation}
and the following identity:
\begin{equation}\labeq{add_der}
\begin{split}
F(\xi , \alpha _1;\tau )
F'_2(\xi , \alpha_2 ;\tau )-
F(\xi , \alpha_2 ;\tau )F'_2(\xi , \alpha_1 ;\tau )\\
=F(\xi , \alpha_1+\alpha_2 ;\tau )&(\wp(\alpha_1)-
\wp(\alpha_2)),
\end{split}
\end{equation}
where $F'_2(\xi , \alpha ;\tau )=\frac{\partial F(\xi , \alpha ;\tau )}{\partial \alpha}$ denotes the derivative with respect to the second argument $\alpha$ and $\wp$ denotes the Weierstra{\ss} function
$\wp =E_2-e_2$. 

In this section about single elliptic curves, $\tau$ will
always be a constant. We will therefore omit it in the notation. 

\sss{Connections and multipliers}
Let $\VC$ be the vector bundle on $X$ defined by $(V,m)$. In the same way as
before, connections on $\VC$ are just the same as 
$\Lambda$-equivariant connections on $\CM$.
 
So let $\nabla$ be a connection on $\VC$. Since the exterior
derivative $d$ on the affine space $\CM$ is itself a $\Lambda$-equivariant
connection, we can write, slightly abusing notation:
$$\nabla = d + \omega,$$
where $\omega$ is a section of the bundle
$\intend{}{\VC}\ot\Omega^1(X)$. 

More explicitely, a connection on $\VC$ is given by its pullback to
$\CM$, which takes the form:
\begin{equation}\labeq{nabla_equiv}
\nabla = d + \omega, \qtext{with}
\omega(\xi+\tau)=m\omega(\xi)m^{-1},\ \omega(\xi+1)=\omega(\xi)
\end{equation}
since $\intend{}{\VC}$ is the bundle with fibre $\endo{}{\VC}$ and 
multiplier the conjugation by $m$. In the general framework of \cite{Andre},
this should be thought as an integrability condition. 

We'll work mostly with multipliers given in exponential form
$m=\exp(n)$. In this setup, the conditions in \refeq{nabla_equiv} become:
\begin{equation}\labeq{adj_eq}
 \omega(\xi+\tau) = \exp(\ad_n)\omega(\xi),\ \omega(\xi+1)=\omega(\xi)
\end{equation}

\begin{defprop}\label{def_conn_fw}
	Let $V$ be a vector space with two
	simultaneously nilpotent operators $\a$ and $A$. The formula
	\begin{equation} \labeq{nabla_fw}
	   \nabla = d - \ad_\a F(\xi,\ad_\a)(A) d\xi,	
	\end{equation}
	where $F$ has to be understood as a Laurent series in the
	second variable,  
	defines a nilpotent meromorphic connection on the vector
	bundle $\VC$ on $X$ with fibre 
	$V$ and multiplier $\exp(-2\pi i\a)$. 
\end{defprop}
\begin{proof} 
We have to check that $\ad_\t F(\xi, \ad_\t)(A)d\xi$ satisfies \refeq{adj_eq} 
for $n=-2\pi i\a$. This is nothing but the 
quasi-periodicity of $F$, as expressed in equations \refeq{F_per} and
\refeq{F_quasiper}. 
\end{proof}

\begin{nota}\label{nota_R}
In the sequel, $R$ will denote the Hopf algebra $\sernc{\CM}{\a,A}$ of formal 
non-commutative series in $\a$ and $A$ over $\CM$, equipped with the 
standard coproduct 
$$
	\Delta(\a)=1\ot\a+\a\ot 1; \Delta(A)=1\ot A + A\ot 1
$$
\end{nota}
 
As a Hopf algebra, $R$ is also the completed universal enveloping algebra of
$\liel{\CM}{\a,A}$, the free Lie algebra generated by $\a$ and $A$.

Vector spaces with simultaneous nilpotent operators $\a$ and $A$ are objects 
of the tensor category $\nmods{R}$ of nilpotent modules over $R$.

We can now formulate the main result of this section:

\begin{thm}\label{thm_equiv_fw}
	The assignment given by the formula \refeq{nabla_fw} 
	extends to an equivalence of rigid tensor categories 
	$$ \FC\colon \nmods{R} \longto \nilpconn(X;e) $$
\end{thm}

Before proceeding to the proof in \ref{sss_plan_thm2} and followings, we'd 
like to expand a bit further on the meaning of this result.
\sss{Tannakian reformulation}
Let us choose once and for all a quasi-inverse of $\FC$ and denote its
composition with the forgetful functor $(V,\a,A)\mapsto V$ by
$\omega_\dr$.
This is a fibre functor, that depends a priori only on the 
choice of the uniformisation $\CM^*\to X$, up to a canonical functorial 
isomorphism.  

The following corollary is then just a restatement of the theorem:
\begin{cor}\label{fw_torsor}
	The fundamental group $\pi_1^\dr(X\setminus(e),\omega_\tau)$ is
	$\exp\liel{}{\a,A}$. 

	The (pro)vector bundle $\PC\ \ass\ \FC(R)$, in which $R$ is
	considered as a left module over itself, 
	is the fundamental torsor $P^\dr(X\setminus(e),\omega_A,-)$. 
\end{cor}
The connection of $\PC$ is given by the formula \refeq{nabla_fw},
in which $\a$ and $A$ have to be interpreted as left multiplications in
$R$. 

\sss{Semi-canonical De Rham paths}
The trivialisation of the pull-back of $\PC$ to $\CM^*$ means moreover that, 
for any $z\in\CM^*\setminus\{1\}$, we have a canonical De Rham path 
between its image $x$ in $X\setminus\{e\}$ and $\omega_A$: that would be the 
$1$ of $\PC_x$; it is simply induced by
the canonical isomorphism between $z^*(V\ot\OC_{\CM^*})$ and $V$. This
applies to tangential base points at $e$ as well. By composition, we get also 
semi-canonical paths between any two points, obeying to similar ambiguities 
upon changes of liftings.

The existence of canonical paths is the
anchoring needed to interpret the parallel transport along a topological path
$\gamma$ for $\PC$ as a De Rham path, and assign a well-defined 
element of $R$ to by comparing it to the canonical path, as in 
\cite[12.15]{DelPi1}. 

These paths really depends on the choice of the lifting: changing $z$ in 
$qz$ amounts to a left multiplication by $\exp(2\pi i\a)$. 

\sss{}
Although $\omega_\dr$ and $\PC$ seem to depend on the choice of
uniformisation, they actually don't. Since this is a special case of the more
general statement for families of section \ref{sec_rel}, we won't repeat
it here. Of course, the explicit description of $\PC$ via its pullback depend 
on this choice as well as of the lifting of base points.  
 
\sss{Plan of the proof}\label{sss_plan_thm2}
The remainder of this subsection is devoted to the proof of theorem
\ref{thm_equiv_fw}. Most of the steps are routine checks and will 
therefore be only sketched. 

We then turn to the essential surjectivity, which just means that 
any connection can be put in the form \refeq{nabla_fw}. To prove this, 
we remark that $A$ can be recovered from the degree $0$ 
component in $\xi$ of $\nabla$, and use proposition 
\ref{prop_unicite_sect} above.

As in the proof of proposition
\ref{thm_bdl}, we reduce  
the question of fully faithfullness to the special case of global sections, 
which is treated first in a separate proposition. 
\sss{Tensor structure}
It is obvious that $\ad_\t F(\xi,\ad_\t)(A)$ is a Lie series
in the variables $\a$ and $A$. Therefore, it is a primitive element with
respect to the Hopf algebra structure of $R$. This ensures that our
functor is compatible with tensor product and duality. 

\sss{Essential surjectivity}\label{sss_fw_ess_surj}
Let $(V,\a)$ represent a vector bundle on $X$ and let $\nabla$ be a
 covariant connection with simple poles at lattice points. We thus have
$\nabla = d - \nu d\xi$ where $\nu$ is a global section of
$\endint{\nilpbdl(X)}{\VC}$ with simple poles at lattice points.
$$ \nu(\xi+\tau)=\exp(-2\pi i\ad_\t)\nu(\xi); \quad \nu(\xi+1)=\nu(\xi) $$

Let's denote by $\nu_0$ the constant term of $\nu$ in its Laurent
series expression at $\xi=0$. 
By proposition \ref{prop_hor_sect}, if we find an endomorphism $A$ of
$V$ such that the degree $0$ term in $\xi$ of $\ad_\t F(\xi,\ad_\t)(A)$ is 
equal to
$\nu_0$, then both sections $\ad_\t F(\xi, \ad_\t)(A)$ and $\nu$ will
have to coincide. 
We have:
\begin{eqnarray*}
   \xi F(\xi,\alpha) &=& 
	\exp(-\sum_{k\geq 1} \frac{(-1)^k\xi^k}{k}(E_k(\alpha)-e_k) \\
&=& 1 + \xi E_1(\alpha) \mod \xi^2, \qtext{hence}\\
   \alpha F(\xi,\alpha) &=& \xi^{-1} + \alpha E_1(\alpha) \mod \xi
\end{eqnarray*}
Therefore, the constant term in $\xi$ 
of $\ad_t F(\xi,\ad_\t)$ is $F_0(\ad_\t)$, 
for $F_0(\alpha)\ \ass\ \alpha E_1(\alpha)$. 
Since $F_0(\alpha) = 1\ \text{mod}\ \alpha$, the operator $F_0(\ad_\t)$ is
invertible, and $F_0^{-1}(\ad_\t)(\nu_0)$ is the $A$ we sought.
\hfill\qed

\begin{prop}[Horizontal sections]\label{prop_hor_sect}
   Let $V$ be an $R$-module. The horizontal sections of $\FC(V)$ are
   the constants annihilated by $\a$ and $A$.	
\end{prop}
\begin{proof}
	Let $s$ be an horizontal section. By proposition
	\ref{prop_hol_sect} above, we already know that $s$ must be
	constant and annihilated by $\a$. Let us express the
	horizontality:
\begin{eqnarray*}
     \nabla s = ds + \alpha F(\xi, \alpha)A\cdot sd\xi &=& 0, \qtext{\ie}\\
  	(\alpha F(\xi, \alpha) A) \cdot s &=& 0, \qtext{hence}\\
	F_0(\alpha)(A) s &=& 0, 
\end{eqnarray*}
   where $F_0$ is as in the proof of essential surjectivity
   above. Since $F_0(\alpha)$ is invertible, we thus obtain the wished $As=0$.
\end{proof}

\sss{Fully faithfullness}\label{fw_full_faith}
Let $V,W$ be $R$-modules. 
An horizontal morphism $f:\FC(V)\to\FC(W)$ is 
an horizontal section of 
$$\inthom{\nilpconn(X;e)}{\FC(V_1)}{\FC(V_2)} = \FC(\inthom{\mods{R}}{V_1}{V_2})$$
By proposition \ref{prop_hor_sect} above, $f$ is a constant linear map
$V\to W$ satisfying $\a f=f\a$ and $Af=fA$, hence a morphism of
$R$-modules. The faithfullness is obvious. \hfill\qed

\sss{The Hodge structure}
The weight and Hodge filtrations on the fundamental De Rham Hopf algebra $R$
can be defined in complete analogy with the genus $0$ case \cite[\S12]{DelPi1}.
\begin{eqnarray*}
     W_n &\ass& \bigoplus_{i\leq -n} R_n\\
     F^p &\ass& \bigoplus_{i\geq -p} \{x\in R, \deg_A(x)= i\},
\end{eqnarray*}
where $R_n$ is the component of degree $n$ of $R$ for the total
degree in $\a$ and $A$. 

It is easy to check that they are compatible with the $\QM$-structure of
\ref{ss_Qsingle}.
 and are invariant with respect to the action of
$\SL_2(\ZM)$ (see \ref{ssect_modul} or definition \ref{def_Rrel}),
hence are well-defined. 

\sss{Simple elliptic polylogarithms}Although we won't detail it here, it is 
possible to show that a regularized version of the generating series 
$\Lambda$ of simple elliptic polylogarithms introduced in \cite{LevinAnal} is 
horizontal for $\nabla_\PC$ in the quotient of $R$ by the second term of its 
derived series. This is in complete analogy with the genus $0$ situation.
 


\let\bfa=\a
\def\cal{\mathcal}

\section{A universal flat connection over the upper-half plane }\label{sec_univ}
\def\t{\mathbf t}
Treating the parameter $\tau$ of the elliptic
curve as a variable, we get the standard family $\Bbb E$ of
elliptic curves over the upper-half-plane
$\Bbb H=\{\tau |{\Im \tau}>0\}$. \joenote{Should I expand on this ?}
Furthermore,
we can define a bundle $\PG$ with fiber $R=\sernc{\CM}{\a,A}$ on $\Bbb E$ by
the same rule as in the case of an individual curve.

In this section, we tie the fiberwise bundles with connection
$(\PC_\tau,\nabla_\tau)$  
from section \ref{sec_fw} together into a \emph{flat} connection
$\nabla_\PG$ on $\PG$ over the \emph{surface} $\Bbb E$. 

\ssect{Combinatorial conventions and lemmas}
We shall make a frequent use of the following special notation: 

\begin{nota}\label{nota_db}
Let $\gG$ a pronilpotent Lie algebra and $t$ in $\gG$. 
For $B,C\in\gG$ and a formal \emph{commutative} power series $f(X,Y)=\sum f_{ij}X^iY^j$, 
we set:
$$ f(X,Y)\db{B,C}_t\ \ass\ \sum_{i,j\geq 0} f_{ij} [\ad_t^iB,\ad_t^jC] $$
\end{nota}
Most of the times, $t$ will be clear from the context, so we'll omit it.

\begin{lemma}[Jacobi identity]\label{lemma_jacobi}
  For any formal power series $f(X,Y)$, we have:
$$\ad_t f(X, Y)\db{B,C}_t = (X+Y) f(X,Y)\db{B,C}_t$$
\end{lemma}

\sss{}\label{sss_db_antisym}
Notice that for a symmetric $f$, i.e., $f(X,Y)=f(Y,X)$, and $C=B$, we have 
$f(X,Y)\db{B,B}=0$, hence in this case one may replace $f$ by its 
antisymmetrisation.

\begin{lemma}\label{lemma_DDB} Let $D$ be a derivation of $R$. For any
  power series $f$ in one variable, we have the following identity:
$$D\, (f(\ad_t)B)= 
\frac{f(X+Y)-f(Y)}{X}
\db{D(t),B}_t +f(\ad_t)D(B)
$$
\end{lemma}
\begin{proof}
Induction on monomials using the Jacobi identity \ref{lemma_jacobi} above.
\joenote{Andrey, do you agree with this cut?}
\end{proof}

Let us also recall this well-known fact \cite[corollary 3.23]{Reut}:
\begin{lemma}\label{lemma_der_exp}
Let $D$ be a derivation of $R$. The following equality holds:
$$D(\exp(t)) = \exp(t)\frac{\exp(-\ad_t)-1}{-\ad_t}D(t).$$ 
\end{lemma}

\ssect{The formula for the connection}
Now we are able to define the connection $\nabla_{\PG}$. As in section 
\ref{sec_fw}, it will be convenient to write Laurent series in $\ad_\t$ within
computations, as long as the polar parts in $\ad_\t$ do eventually vanish.
\begin{eqnarray*}
\psi_A(\tau) &\ass& -\frac12\frac{1}{2\pi i}
\frac{XY}{X+Y}
\left(\left(\wp(X,\tau)-\frac1{X^2}\right)-\left(\wp(Y,\tau)-\frac1{Y^2}\right)
\right)
\db{A,A}d\tau\\
\psi_{\t}(\tau) &\ass& -\frac1{2\pi i} Ad\tau,\\
\nu(\xi,\tau) &\ass& 
-\left(\ad_\t F(\xi , \ad_\t ;\tau ) d\,\xi
+\frac{1}{2\pi i}\left(\ad_\t F'_2(\xi , \ad_\t ;\tau )+
\frac{1}{\ad_\t}\right)d\,\tau\right)A,
\end{eqnarray*}
and, finally, 
\begin{equation}\labeq{def_P_rel}
\nabla_{\PG}\ \ass\ d  +\nu  
+\psi_A 
\frac{\partial}{\partial A}+\psi_{\t}
\frac{\partial}{\partial\t},
\end{equation}
where
$f\frac{\partial\,}{\partial\, A}$ (resp. 
$f\frac{\partial\,}{\partial\, \t }$) has to be understood as the
derivation which maps $A$ to $f$ and $\a$ to $0$ (resp. $A$ to $0$ and
$\a$ to $f$). \joenote{Andrey, in the context of the
Grothendieck-Teichmüller group for genus $0$, similar operators
  are called special derivations. This notation has the advantage of
  being light, but it's also quite
  confusing. We should maybe find something else.}
Hence the commutator 
$[f\frac{\partial\,}{\partial\, A},r ]$
 of $f\frac{\partial\,}{\partial\, A}$
and $r $ is just $f\frac{\partial\,r}{\partial\, A}  $

\begin{prop}\label{univ_nabla_quasiper}
 The operator $\nabla_{\PG}$ above descends to  a connection 
 on the vector bundle $\PG$ on $\EM$. 
\end{prop} 
\begin{proof}
We must check that $\nabla_{\PG}$
is invariant with respect to the shift $\xi\to\xi +1$
and that it transforms in correct
way under $\xi \to \xi+\tau$. The first 
property
follows immediately from the invariance of $F$ with
respect to this shift. So we shall check the quasiperiodicity with respect 
to $\tau$:
\begin{equation}\labeq{univ_P_twist}
\nabla_{\PG}|_{\xi+\tau}=\exp(-2\pi i \t )\nabla_{\PG}|_{\xi }
\exp(2\pi i \t ).
\end{equation}
The $\nu$ summand in the left hand side of \refeq{univ_P_twist} is equal to  
$$
\begin{array}{rl}
-&\left(\ad_\t F(\xi+\tau , \ad_\t ;\tau ) d\,(\xi+\tau)
+\frac{1}{2\pi i}\left(\ad_\t F'_2(\xi +\tau, \ad_\t ;\tau )+
\frac{1}{\ad_\t}\right)d\,\tau\right)A \\
-&\left(\exp(-2\pi i\ad_\t)(\ad_\t F(\xi  , \ad_\t ;\tau ) d\,\xi 
 + \ad_\t F(\xi  , \ad_\t ;\tau )d\,\tau) 
\phantom{\frac{1}{\ad_\t}}\right.\\
+& \left.\frac{1}{2\pi i}\left(\exp(-2\pi i\ad_\t)(
-2\pi i\ad_\t F(\xi  , \ad_\t ;\tau )+\ad_\t F'_2(\xi, \ad_\t ;\tau ))+
\frac{1}{\ad_\t}\right)d\,\tau\right)A\\
=&-\exp(-2\pi i\ad_\t)\left(\Aug\ad_\t F(\xi  , \ad_\t ;\tau ) d\,\xi \right.\\
&\left.\hspace{7.5em}
 + \left( \frac{1}{2\pi i}\left(
 \ad_\t F'_2(\xi, \ad_\t ;\tau )+
\frac{1}{\ad_\t}\right)+
\frac{\exp(2\pi i\ad_\t)-1}{2\pi i\ad_\t}\right)
d\,\tau\Aug\right)A
\end{array}$$
The other terms in the left hand side of \refeq{univ_P_twist} 
do not depend on $\xi$; they are therefore invariant.
In the right hand side, the term $\exp(-2\pi i \t )\nu \exp(2\pi i \t )$ 
is equal to
$$ -\exp(-2\pi i\ad_\t)\left(\ad_\t F(\xi  , \ad_\t ;\tau ) d\,\xi 
 +   \frac{1}{2\pi i}\left(
 \ad_\t F'_2(\xi +\tau, \ad_\t ;\tau )+
\frac{1}{\ad_\t}\right)  
d\,\tau\right)A,
$$
whereas the term 
$\exp(-2\pi i  \t ) \psi_A \frac{\partial}{\partial A} 
\exp(2\pi i \t)$ equals $\psi_A \frac{\partial}{\partial A}$ since 
$\psi_A\partial_A(\tau)=0$.

The last remaining term 
$\exp(-2\pi i\t) \psi_{\t}  \frac{\partial\,}{\partial\,  \t } 
\exp(2\pi i  \t )$ is equal to
$$\exp(-2\pi i\t)\left(-\frac1{2\pi i}A \frac{\partial}{\partial\t}\right) 
\exp(2\pi i\t) = \frac{\exp(-2\pi i\ad_\t)-1}{-2\pi i\ad_\t} 
(- A)-\frac1{2\pi i}A \frac{\partial}{\partial\t},
$$
by lemma \ref{lemma_der_exp}. So, both sides of \refeq{univ_P_twist} are indeed
equal.
\end{proof}
\begin{prop}
The connection $\nabla_{\PG}$ is flat.
\end{prop}
\begin{proof}
Denote by $\omega$ the differential form
$\nu  
+\psi_A 
\frac{\partial}{\partial A}+\psi_{\t}\frac{\partial}{\partial\t}.$
 
We shall prove that $d\,\omega+\omega\wedge\omega=0$.
The term $d\omega$ is equal to $d\nu $, as other terms are lifted from
the one dimensional $\Bbb H$, and therefore are closed. 
\begin{eqnarray*}
d\omega&=&-d\, \left(\ad_\t F(\xi , \ad_\t ;\tau ) d\xi
+\frac{1}{2\pi i}\left(\ad_\t F'_2(\xi , \ad_\t ;\tau )+
\frac{1}{\ad_\t}\right)d\,\tau\right)A \\
&=&\left(\frac{\partial}{\partial\tau}
         (\ad_\t F(\xi , \ad_\t ;\tau )) 
         -\frac{1}{2\pi i}\frac{\partial}{\partial\xi}
             \left(\ad_\t F'_2(\xi, \ad_\t; \tau )+\frac{1}{\ad_\t}\right)
    \right)Ad\xi\wedge d\tau \\
&=&
\ad_\t\left( 
\frac{\partial}{\partial\,\tau}
F(\xi , \ad_\t ;\tau ) -\frac{1}{2\pi i}\frac{\partial}{\partial\,\xi}   F'_2(\xi , \ad_\t ;\tau )
  \right)A d\xi\wedge d\tau = 0,
\end{eqnarray*}
as the function $F$ satisfies the mixed heat 
equation \refeq{mixed_heat}.
 
Let's now turn our attention to the $\omega\wedge\omega$ term. We have
\begin{eqnarray*}
\omega\wedge\omega &=& 
(2\pi i)^{-1} (\Sigma_1+\Sigma_2+\Sigma_3)d\xi\wedge d\tau, \qtext{with}\\
\Sigma_1&\ass&\left[\ad_\t F(\xi , \ad_\t ;\tau ) A,
\left(\ad_\t F'_2(\xi , \ad_\t ;\tau )+
\frac{1}{\ad_\t}\right)A\right] \\
\Sigma_2 &\ass& -\left[\frac12\frac{XY}{X+Y}
    \left(\left(\wp(X)-\frac1{X^2}\right)
    -\left(\wp(Y)-\frac1{Y^2}\right)
    \right)\db{A,A}
    \frac{\partial}{\partial\, A},\right. \\
&&\left.\Aug\qquad \ad_\t F(\xi , \ad_\t ;\tau ) A \right] \\
\Sigma_3&\ass& \left[ A \frac{\partial}{\partial\t},
           \ad_\t F(\xi , \ad_\t ;\tau ) A \right] 
\end{eqnarray*}
Till the end of this computation, it will be safe to omit $\tau$ from the 
notation. Using notation \ref{nota_db} and remark \ref{sss_db_antisym} on the 
first summand, we have: 
\begin{gather*}
\Sigma_1 = XF(\xi , X) \left(Y F'_2(\xi , Y)+ \frac{1}{Y}\right)\db{A,A}\\
= \frac12\left(XF(\xi , X)
    \left(Y F'_2(\xi , Y ;\tau )+ \frac{1}{Y}\right)
          -Y F(\xi , Y) \left(XF'_2(\xi , X)+\frac{1}{X}\right)
    \right)\db{A,A}  
\end{gather*}
Applying Jacobi identity in the form of lemma \ref{lemma_jacobi} to $\Sigma_2$, we obtain:
$$
\Sigma_2 = 
   -\frac12(X+Y) F(\xi, X+Y) \frac{XY}{X+Y}
   \left(\left(\wp(X)-\frac1{X^2}\right)
         -\left(\wp(Y)-\frac1{Y^2}\right)\right)\db{A,A}
$$
For the third summand, we apply lemma \ref{lemma_DDB} and get:
\begin{eqnarray*}
\Sigma_3 &=& -\frac{(X+Y) F(\xi , X+Y)- Y F(\xi ,  Y)}{X} 
\db{A,A} \\
&=& -\frac12\left(\frac{(Y^2-X^2) 
F(\xi , X+Y)+X^2 F(\xi ,  X)-Y^2 F(\xi ,  Y)}{XY}\right)
\db{A,A}
\end{eqnarray*}
Putting the three summands together again, we have:
\begin{eqnarray*}
\omega\wedge\omega &=& 
2^{-1}(2\pi i)^{-1} \Xi(X,Y,\tau)\db{A,A}d\xi\wedge d\tau, \qtext{with}\\
\Xi(X,Y,\tau) &\ass& XF(\xi , X ) 
 \left(Y F'_2(\xi , Y)+ \frac{1}{Y}\right)
-Y F(\xi , Y) \left(X F'_2(\xi , X)+ \frac{1}{X}\right) \\
&& - {XY}   F(\xi , X+Y)  \left(\left(\wp(X)-\frac1{X^2}\right)
-\left(\wp(Y)-\frac1{Y^2}\right)
\right) \\
&&- \frac{(Y^2-X^2) F(\xi , X+Y)+ X^2 F(\xi ,  X)
          - Y^2 F(\xi ,  Y)}
{XY}\\
&=& XF(\xi , X) Y F'_2(\xi , Y) -Y F(\xi , Y)  X
F'_2(\xi , X) \\
&&- XY F(\xi , X+Y)   \left(\wp(X) -\wp(Y)\right).
\end{eqnarray*}
This latter expression vanishes, according to \refeq{add_der}; so does the curvature of $\nabla_{\PG}$.
\end{proof}
\ssect{Modularity}\label{ssect_modul}
\sss{}
Consider the standard action of the modular group 
$SL_2(\Bbb Z)$ on the upper-half-plane $\Bbb H$ and
the relative curve $\Bbb E$ over $\Bbb H$:
$\tau\to \tau' =\frac{a\tau +b}{c\tau +d};\quad
\xi\to \xi'=\frac{\xi}{c\tau +d}$. 
Denote by ${\cal R} $ the (pro)bundle on $\Bbb H$ with fiber $R$.
Define an action of $SL_2(\Bbb Z)$ on $\cal R$
by the rule:
$\t \to \t'=\frac{\t }{c\tau +d}$,  
$A\to A'=(c\tau +d)A +2\pi i c \t$; and  lift this action to an action on $\PG$ by  multiplier 
$\exp(2\pi i\frac{c \xi \t}{c\tau +d})$. 
Thus, a $SL_2(\Bbb Z)$-invariant section of $\PG$   is
an $R$-valued function $f$ in $\tau$ and $\xi$
such that 

\begin{equation}f\left(\frac{a\tau +b}{c\tau +d}, \frac{\xi}{c\tau +d};
\frac{\t }{c\tau +d}, (c\tau +d)A +2\pi i c\t)=
\exp\left(2\pi i\frac{c \xi \t}{c\tau +d}\right) f(\tau, \xi;
\t,A\right).
\end{equation}

This action is well defined (compatible with the group structure on $SL_2(\Bbb Z)$ and compatible
with fiberwise action of $\Bbb Z^2$).
We leave these as an exercise to the reader. 

As before, the condition that a global endomorphism has to fulfill is obtained 
by replacing $\t$ with $\ad_\t$.

\begin{prop}The connection $\nabla_{\PG}$
is equivariant for the action of $SL_2(\Bbb Z)$ on the vector bundle $\PG$.
\end{prop} 
\begin{proof}
$$\nabla_{\PG}  |_{\tau', \xi';
\t',A'}
=\exp\left(2\pi i\frac{c \xi \t}{c\tau +d}\right)\nabla_{\PG} |_{\tau, \xi;
\t,A}\exp\left(-2\pi i\frac{c \xi \t}{c\tau +d}\right).$$

The term $\nu$ in the l.h.s. equals
$$
-\left(\frac{\ad_\t}{c\tau+d} 
F\left(\frac{\xi}{c\tau+d} , 
\frac{\ad_\t}{c\tau+d} ;\frac{a\tau +b}{c\tau+d} \right) 
d\,\frac{\xi}{c\tau+d}+\right.$$
$$\left.
\frac{1}{2\pi i}\left(
\frac{\ad_\t}{c\tau+d} 
F'_2(\frac{\xi}{c\tau+d} , 
\frac{\ad_\t}{c\tau+d} ;\frac{a\tau +b}{c\tau+d} 
 )+
\frac{{c\tau+d}}{ {\ad_\t}}\right)
d\,\frac{a\tau +b}{c\tau+d} \right)((c\tau +d)A+
2\pi i c\t)=
$$
 $$
-\exp\left(2\pi i\frac{c \xi \ad_\t}{c\tau +d}\right)\left( {\ad_\t} 
F( {\xi}  , 
 \ad_\t  ; \tau   ) 
\left(d\, {\xi}-\frac{c\xi}{c\tau+d}d\,\tau\right) -2\pi ic\t\left(\frac{d\, \xi}{c\tau+d}-\frac{c\xi 
d\,\tau}{(c\tau+d)^2}\right) +\right.$$
$$\left.
\frac{1}{2\pi i}\left(
 \ad_\t \left(2\pi i\frac{c\xi}{c\tau+d}F( {\xi}  , 
 \ad_\t  ; \tau   ) +  
F_2'( {\xi}  , 
 \ad_\t  ; \tau   ) \right)+
\exp\left(-2\pi i\frac{c \xi \ad_\t}{c\tau +d}\right)\frac{1}{ {\ad_\t}}\right)
d\,\tau \right)A
 =
$$
$$\exp\left(2\pi i\frac{c \xi \ad_\t}{c\tau 
+d}\right)
\nu +\frac{\exp\left(2\pi i\frac{c \xi \ad_\t}{c\tau +d}\right)-1}{\ad_\t}d\,\tau A-2\pi ic\t\left(\frac{d\, \xi}{c\tau+d}-\frac{c\xi
d\,\tau}{(c\tau+d)^2}\right).$$

In this calculations we use the modular property
of the function $F$ and the evident equality
$f(\ad_\t)\t=f(0)\t$.

We have  $$\frac{\partial}{\partial A'}=
\frac1{c\tau +d} \frac{\partial}{\partial A},\quad\frac{\partial}{\partial \t'}=
(c\tau+d)\frac{\partial}{\partial \t }
-2\pi i c\frac{\partial}{\partial A}.$$ 
 As $\psi_A$ is a modular form of weight $1$:
$$\psi_A(\tau', \xi'; \t',A') 
=-\frac12\frac{1}{2\pi i}(c\tau +d)
\frac{XY}{X+Y}
\left(\left(\wp(X)-\frac1{X^2}\right)
-\left(\wp(Y)-\frac1{Y^2}\right)
\right)\times$$
$$
((c\tau +d)A +2\pi i c \t,(c\tau +d)A +2\pi i c \t)\frac{d\,\tau}{(c\tau+d)^2} =(c\tau+d)
\psi_A(\tau, \xi;
\t,A ),$$
the term $\psi_A\frac{\partial}{\partial A}$ is $SL_2(\Bbb Z)$ invariant.

We have also:
$$\psi_{\t}(\tau', \xi';
\t',A') =-\frac1{2\pi i}((c\tau +d)A +2\pi i c \t)\frac{d\,\tau}{(c\tau+d)^2}=\frac1{c\tau +d}
\psi_{\t} -  c \t\frac{d\,\tau}{(c\tau+d)^2}.
$$

Note that $A'$ and $\t'$ depend in $\tau$ 
explicitly as function in $A$ and $\t$, so
$$ \left.\frac{ d}{d\,\tau}\right|_{ \t',A'}=  \frac{\partial }{\partial\,\tau} -\frac{cA'}{(c\tau +d)^2}
\frac{\partial }{\partial\,A}
+c\t'\frac{\partial }{\partial\,\t}
$$
Therefore
$$d|_{ \t',A'} +\psi_{\t}  
({\tau', \xi';
\t',A'}) 
\frac{\partial\,}{\partial\, \t' }= d|_{\t, A}+ 
d\tau\left(
  -\frac{cA'}{(c\tau +d)^2}
\frac{\partial }{\partial\,A}
+c\t'\frac{\partial }{\partial\,\t}\right)$$
$$\left(\frac1{c\tau +d}
\psi_{\t} -  c \t\frac{d\,\tau}{(c\tau+d)^2} \right)
(c\tau+d)\frac{\partial}{\partial \t }+$$
$$\left(-\frac1{2\pi i}A'\frac{d\,\tau}{(c\tau+d)^2}
\right)\left(-2\pi i c\frac{\partial}{\partial A}\right)=d|_{\t, A} +\psi_{\t}
\frac{\partial}{\partial \t }.
$$

In the r.h.s. we have:
$$\exp\left(2\pi i\frac{c \xi \t}{c\tau +d}\right)\left(d -\frac1{2\pi i} Ad\,\tau\frac{\partial\,}{\partial\, \t }
\right)\exp\left(-2\pi i\frac{c \xi \t}{c\tau +d}\right)=$$
$$d  -
2\pi ic\t\left(\frac{d\, \xi}{c\tau+d}-\frac{c\xi 
d\,\tau}{(c\tau+d)^2}\right) +\frac{\exp\left(2\pi i\frac{c \xi \ad_\t}{c\tau +d}\right)-1}{\ad_\t}d\,\tau A,
$$
the second summand comes from the differentiation 
by $\xi$ and $\tau$, the last can be calculated
using Lemma \ref{lemma_DDB}.  

Evidently,
$$
\exp\left(2\pi i\frac{c \xi \t}{c\tau +d}\right)\nu\exp\left(-2\pi i\frac{c \xi \t}{c\tau +d}\right)=\exp\left(2\pi i\frac{c \xi \ad_\t}{c\tau 
+d}\right)
\nu;$$
and $\psi_A 
\frac{\partial\,}{\partial\, A}$ commutes with
$\exp\left(-2\pi i
\frac{c \xi \t}{c\tau +d}\right)$.
This finishes the proof.
\end{proof}

Note that such a lifting of the bunch of fiberwise connections to a 
flat connection is not unique. We shall return to this topic in \ref{ss_geom}.

Indeed, one can
twist the connection $\nabla_{\PG}$ by any $R$-valued $SL_{2}(\Bbb Z)$-invariant form  $\mu$ on $\Bbb H$
by the right multiplication on $\mu$:
$$\nabla_{\PG}^{\mu}(G)=\nabla_{\PG}(G)+G\mu=\nabla_{\PG}(G)+\mu G- {\rm ad}_\mu G.$$
This connection is flat as the right and the left multiplications commute.

\section{Equivalence for the relative case}\label{sec_rel} 
	In this section, we generalize the results of the previous section to
the case of a smooth family $p: X\to S$ of elliptic curves.  We refer to the
case treated in the previous section as the ``fiberwise case''. We
denote by $e$ the unit section of $X$. 

We are mostly interested in the following families:
\begin{enumerate}
	\item $S$ is Poincaré's upper-half plane $\HM$ with coordinate
	$\tau$ and $X$ is the standard family $\EM$.
	\item $S$ is the punctured unit disc $D^*$ with coordinate $q$ and
	$X$ is $(\CM^*\times D)/q^\ZM$.  
	\item $S$ is some modular curve, or even the moduli stack
	$\M{1,1}$ of elliptic curves (smooth genus $1$ curves with one
	marked point) 
\end{enumerate} 

We provide here an explicit description of the \emph{relative} nilpotent
de Rham fundamental group and torsor of $X$. These objects are defined
in \ref{ss_rel_def}. A precise statement is given as theorem \ref{thm_rel}. 

\ssect{Generalities}\label{ss_rel_def}

\begin{defi}
	A connection will be said to be relatively trivial on $X$ if
	it is a pullback from $S$ by $p$. 

	A relatively nilpotent connection on $X$ (resp., on
	$X\setminus e$, on $X$ with simple poles at $e$) is an
	iterated extension of relatively trivial connections, in the
	relevant category of flat connections. 
\end{defi}

We denote by $\nilpconn(X/S;e)$ the category of relatively
nilpotent flat meromorphic connections on $X$ wich are holomorphic
on $X\setminus e(S)$ and have simple poles at the unit section, 
and by $\conn(S)$
the category of \emph{all} flat connections on $S$. These are tannakian 
categories, and the pullback functor
$$ p^*\colon \conn(S) \longto \nilpconn(X/S; e) $$
respects the tannakian structure. 
\sss{Relative tannakian theory}\label{sss_rel_tann}
We give here a short account of some results of Spitzweck \cite{spitzdiplom},
generalizing the theory of tannakian categories to such
relative situations. 

In general, one considers a
tensor functor $i\colon \SC\to\TC$ between two tannakian
categories. \emph{Relative fiber functors are} defined similarly as the usual ones,
except that they take values in categories of vector bundles over
$\SC$-schemes, which are also defined in loc.cit., and have to be 
compatible with $i$. 

The category $\TC$ is said to be \emph{neutral} relative to $i$ if there exists a
relative fixed fiber functor over $\SC$. Such an object boils down to a 
pair $(\omega, \phi)$ where $\omega$ is a tensor functor from $\TC$ to
$\SC$ and $\phi:\omega i\isomto\id$. In this case, there is an
$\SC$-affine group scheme $\pi_1(\TC/\SC, \omega)$ and an equivalence
of categories:
\begin{equation}
   \TC \isomto \mods{\pi_1(\TC/\SC, \omega)}
\end{equation}
Further, there exists a pro-coalgebra in $\TC$, denoted by $\PC(\TC/\SC, \omega,
-)$, called the 
\emph{fundamental torsor coalgebra relative to $i$ at $\omega$}. It comes equipped
with a right action of $i\pi_1(\TC/\SC, \omega)$,  and is characterized by
the fact that the above equivalence is given by the
twisting:\joenote{There, I should care about what depends on the
unipotence assumption.}
\begin{equation}
   T \longmapsto \PC(\TC/\SC, \omega, -)\ot_{i\pi_1(TC/SC,\omega)} i(T),
\end{equation}
exactly as in the usual case. 

As a side remark, let us observe that any tannakian category $\SC$
over a field $k$ is relative over $\vects{k}$. There is indeed \cite[2.9]{DelFest}
a canonical tensor functor $\vects{k}\to\SC$. In geometric categories
such as $\conn(S)$, it is given by trivial objects, hence is simply the
pullback functor.  

\begin{defi}\label{def_PG}
For a fixed relative fiber functor $\omega$ on $\nilpconn(X/S;e)$ with
values in $\conn(S)$, 
we will denote by $\pi_1^\dr(U/S, \omega)$ and 
$\PG^\dr(U/S,\omega)$ the fundamental group and torsor of
$\nilpconn(X/S,e)$ at $\omega$, relative to the pullback functor $p^*$. 
\end{defi}
These objects are respectively a pro-Hopf algebra in $\conn(S)$ and a 
pro-coalgebra objects in $\nilpconn(X/s;e)$:  flat connections equipped with 
structure laws which are expressed by diagrams of
horizontal morphisms.  \joenote{Here, I try to recatch the
attention of people not so much delighted with abstract
nonsense. Really useful or brings confusion ?}

\sss{Remarks}
A typical example of a relative fiber functor over $\conn(S)$ would be
the pullback functor by a section of $p:X\setminus\{e\}\to S$. 
It is therefore in general not a priori granted 
that relative fiber functors exist, i.e. that $\nilpconn(X/S;e)$ is
neutral over $\conn(S)$. 

On the other hand, for an {\'e}tale covering $S'\to S$, any section of
$X_{S'}$ provides a fixed fiber functor over $S'$.   
A relative tannakian theory restricted to such geometric fiber functors
has been formulated and used by Wildeshaus \cite[I.3]{WildesMixedFal}.   

\subsection{The main statement}
\sss{A connection algebra on $\M{1,1}$}
In this subsection, we define for any $S$ a Hopf algebra in $\conn(S)$
which will turn out to be the fundamental Lie algebra of
$\nilpconn(X/S;e)$. We start with the case of the upper-half plane.
\begin{defi} \label{def_Rrel}
Let $\RC_\HM$ be the trivial (pro)bundle on $\HM$ of Hopf algebras with fiber 
  $R=\sernc{\CM}{\a,A}$, equipped with the connection
$$\nabla_\RC\ \ass\ d 
+\psi_A \d_A + \psi_\a \d_\a,$$ where we put, using notation \ref{nota_db}:
\begin{eqnarray}  
\psi_A &\ass& 
-\frac12\frac{1}{2\pi i}
\frac{XY}{X+Y}
\left(\left(\wp(X)-\frac1{X^2}\right)
-\left(\wp(Y)-\frac1{Y^2}\right)
\right)\odb A,A\cdb d\tau,\\
\psi_{\t}&\ass&-\frac1{2\pi i} Ad\tau,
\end{eqnarray}
This bundle comes further equiped with the following action of $\SL_2(\Bbb Z)$:
\begin{eqnarray}
 \t \to \t' &=& \frac{\t }{c\tau +d}  \\
 A\to A'&=& (c\tau +d)A +2\pi i c \t
\end{eqnarray}
\end{defi}
This $\RC_\HM$ is actually a Hopf algebra in the category of 
connections on $\HM$:  
\begin{prop} 
	The connection $\nabla_\RC$ is compatible with the Hopf 
	algebra structure of $\RC_\HM$ and $\SL_2(\ZM)$ equivariant.  
\end{prop} 
\begin{proof}
	This follows directly from the fact that $\psi_\a\d_\a +
	\psi_A\d_A$ is a derivation which maps $\a$ and $A$ to
	primitive elements and the freeness of $R$ as an
	algebra.

	The $\SL_2(\ZM)$ equivariance is immediate.
\end{proof} 
\sss{Local uniformisations} 
Let $p:X\longto S$ be a smooth family of elliptic curves. 
Any local symplectic basis of $R^1p_*(\ZM_X)$
over an open 1-connected subvariety $U$ of $S$ defines a morphism $\tau: U\to\HM$,
and an isomorphism $X_{|U}\isomto \tau^*\EM$, where $\EM$ denotes the
standard family over $\HM$. We call such a pair $(U,\tau)$ a local
uniformisation of $X$.\joenote{So the local isomorphism with the
pullback is implicit} 
\begin{prop}\label{prop_Rrel}
	Let $X\longto S$ be a familly of elliptic curves. There is a
	well-defined Hopf algebra $\RC_S$ in $\conn(S)$ such that, for
	any local uniformisation ($U$, $\tau$): 
	\begin{equation}
	 \RC_S|U \isomto  \tau^*\RC
	\end{equation}
\end{prop}
\begin{proof}[Sketch of proof]
	The $\SL_2(\ZM)$ invariance of $\RC$ precisely provides a
	canonical isomorphism between the pullbacks by two local
	uniformisations. The cocycle condition is nothing but the
	associativity of the action. 
\end{proof}
Of course, this is nothing but a reassertion of the well-known fact
that the analytic stack $\M{1,1}$ is the quotient of $\HM$ by
$\SL_2(\ZM)$ in the sense of stacks, together with the description of
sheaf-like objects on quotient stacks by equivariance.

\sss{Remarks} Note that the term  $\psi_{\t}
\frac{\partial}{\partial\t}$ corresponds to the connection on the
tensor algebra of the first homology group equipped with Gauß-Manin
connection. 

One can also consider the (sub)bundle $\cal L$ of free Lie algebras generated 
by $\t$ and $A$. It is invariant under $SL_2(\Bbb Z)$ and the connection 
$\nabla_{\RC}$.

\begin{prop}  
On any family $p:X\to S$ of elliptic curves,  there is a well-defined
relatively nilpotent connection $\PG_S$, such that, for any local
uniformisation 
$(U,\tau)$, 
$$ \PG_S\isomto\tau^*\PG_\HM. $$
Moreover, $\PG_S$ comes equipped with a right $p^*\RC_S$ module
structure, which is locally the pullback of the right
$p^*\RC_\HM$-module structure on $\PG_\HM$.
\end{prop}
\begin{proof} 
	As in proposition \ref{prop_Rrel}, this is a simple
	consequence of the $\SL_2(\ZM)$ equivariance of $\PG$ and of
	the action morphism $\PG\ot p^*\RC\to\PG$. 
\end{proof}

We can now state the main result of this section.
\begin{thm}\label{thm_rel}
	For any smooth family of elliptic curves $p:X\to S$, the functor 
\begin{equation}\labeq{def_F}
\FC_S\colon\left.\begin{array}{lcl}\mods{\RC}&\longto&\nilpconn(X/S;e)
\\
\VC&\longmapsto&\PG\ot\limits_{p^*\RC}p^*\VC \end{array}\right.
\end{equation}	
is an equivalence of categories. 
\end{thm}
One may find the following reformulation in the relative tannakian terms
of \ref{sss_rel_tann} to be more telling:
\begin{cor}
	there is a well-defined and fixed fiber functor  
	$$\omega_S\colon\nilpconn(X/S;e)\to\conn(S),$$
	which makes $\nilpconn(X/S;e)$ neutral relative to $p^*$, and such that: 
	\begin{enumerate}
	\item the fundamental Hopf algebra at $\omega_S$ is $\RC_S$.	
	\item the relative fundamental torsor
	$P^\dr((X\setminus e)/S, \omega_S)$ is $\PG_S$. 
	\end{enumerate}
\end{cor}
As before, $\omega_S$ is the composite of an inverse of $\FC_S$ and
the forgetful functor $\mods{\RC}\to\conn(S)$. 

\sss{Explicit form} \label{sss_expl_rel}
The pullback condition for a local uniformisation
 	$(U,\tau)$ provides directly an explicit description of
 	$\PG$ on $U$. Since the formulas involving the function $F$
 	can all be written in terms of $q$, things stay explicit under
 	the milder assumption 
 	that $X\to S$ admits a Schottky
 	uniformisation, i.e., a function $q:S\to D^*$, such
 	that $X\to S$ is the pullback of $*\CM^*\times S)/q^\ZM$ by
 	$q$, which amounts to the choice of a section of
 	$R^1p_*(\ZM_X)$.  

	Indeed, to a pair $(\VC,\a)$, where $\VC$ is
a vector bundle on $S$ and $\a$ 
is an endomorphism of $\VC$, we can, as in \ref{ss_bdl_fw}, 
associate a nilpotent bundle $\GC(\VC,\a)$ on $X$:
\begin{equation}\labeq{glob_multiplier}
\Gamma(U, \GC(\VC, \a))\ \ass\ \{s\in \Gamma(\pr_2^{-1}U,
\pr_2^*\VC), s(qz) = \exp(-2\pi i \a) s(z)\}, 
\end{equation}
where $\pr_2$ is the projection to the second factor of $\CM^*\times S$. 
More precisely, we have a functor 
$$\GC: (\text{vector bundles with nilpotent endomorphism on $S$})\longto\nilpbdl(X/S),$$
As before, we'll refer to $\GC(\VC,\a)$ as \emph{the bundle with
multiplier $\exp(-2\pi i\a$}). 

	The fundamental torsor $\PG_S$ can then be described as the bundle with
	fiber $R$ and multiplier 
	$\exp(-2\pi i\a)$, 
        equipped with the connection whose pullback to $\CM^*\times S$ is:
	\begin{eqnarray}\labeq{nablaP}
	   \nabla &\ass& d - \nu + \psi_\a\d_\a + \psi_A\d_A, 
	\qtext{where}\\
	\nu  &\ass&  \ad_\a F(q , \ad_\a ;q) \frac{dz}{z}
	+\frac{1}{2\pi i}\left(\ad_\t F'_2(q , \ad_\a ;q)+
	\frac{1}{\ad_\t}\right)\frac{dq}{q}A, 
	\end{eqnarray}
	and $\psi_\a,\psi_A$ are as in definition
	\ref{def_Rrel}. 	

\sss{Reformulation for moduli stacks}
	The theorem can be summarized in the following way: 
	
	There is a canonical fixed fiber
	functor $\omega$ on $\nilpconn((\X{1,1}\setminus
	e)$, relative to $\conn(\M{1,1})$, where $\X{1,1}$ is the
	universal family of 
	elliptic curves. 

	 The relative fundamental Hopf algebra
	$\U\pi_1^\dr((\X{1,1}\setminus e)/\M{1,1}, \omega)$ is $\RC$ and
	the fundamental torsor is $\PG$. Note also that
	$\X{1,1}\setminus e$ can be identified with $\M{1,2}$, the
	moduli stack of smooth curves of genus $1$ with two marked points. 

\ssect{Proof of theorem \ref{thm_rel}}
We treat the local situation in \ref{local_ess_surj1} to \ref{local_ess_surj2}.
We explain how to glue those equivalences in \ref{sss_rel_glue}, 
which is nothing but the well-known 
statement that isomorphicity for stacks can be checked locally. 

\sss{Essential surjectivity in the local case}\label{local_ess_surj1}

In this part, we work over a local uniformisation $(U,\tau)$, 
where $U$ is contractible. We'll provide full details for the 1-dimensional 
case, i.e, essentially $U=\HM$. 

Let $(\VG, \nabla)$ be a object of $\nilpconn(X/U;e)$. According to a 
relative version of Atiyah theorem, the category of vector bundles on $X$, 
relatively  nilpotent over $U$ is equivalent to the category of modules over 
$\serc{\OC_U}{\t}$, and we have a description of the underlying bundle of 
$\VG$ as a trivial bundle $\VC$ on $U$ with fiber $V$ and 
multiplier $exp(\t(\tau))$ in the $\xi$ direction. 

By using the same formal inverting procedure as in \ref{sss_fw_ess_surj}, 
we get a section $A(\tau)$ such that each fiber of $\FC_U(\VC, \t, A)$ 
coincides with our $(\VG, \nabla)$. In other words, its pull-back to 
$\CM^*\times U$ can be written as:
\begin{equation}\labeq{nabla_rel_ess_surj}
   \nabla = d - \nu(\t, A) - \phi(\xi, \tau)d\tau
\end{equation}
Note that $\phi$ must be holomorphic, since the poles of $\nabla$ have 
$\dlog\xi$ form.

\begin{prop}\label{prop_phi_cst}
The expression $\phi(\xi, \tau)$ does not depends in $\xi$.\\
We have $\frac {d\t}{d\tau} - [\phi, \t] = -\frac A{2\pi i}$.
\end{prop}
\begin{proof}
Writing down the quasiperiodicity with respect to $\xi\to \xi+\tau$ of 
$\nabla$ and taking into account the computations from 
\ref{univ_nabla_quasiper}, we get:
\begin{equation}\labeq{phi_quasiper}
\phi(\xi+\tau, \tau)-\exp(-2\pi i\ad_\t)\phi(\xi, \tau)=
-\frac{\exp(-2\pi i\ad_\t)-1}{-\ad_\t}
\left(\frac d{d\,\tau}\t +\frac A{2\pi i}\right)
\end{equation}

Denoting for a while by $K$ the expression 
$\left(\frac{d\t}{d\tau}+\frac A{2\pi i}\right)$, let us observe that function 
$-(F(\xi, \ad_\t)-\frac 1\ad_\t)K$ also satisfies to \refeq{phi_quasiper}, 
so
$\phi(\xi,\tau)+(F(\xi, \ad_\t)-\frac 1\ad_\t)K$ satisfies to the 
corresponding homogeneous equation. By \ref{sss_fw_ess_surj}\joenote{Single out as a lemma or prop?}, it can therefore be 
written as $\ad_\t F(\xi, \ad_\t)C(\tau)$ for some $C(\tau)$. This gives us:
\begin{equation}\labeq{phi_cst}
 \phi(\xi, \tau)=\ad_\t F(\xi, \ad_\t)C(\tau)-(F(\xi, \ad_\t)-\frac 1\ad_\t)K
\end{equation}
Since $\phi$ is holomorphic, the singular part at $\xi=0$ of the latter 
expression must vanish. This provides first $K=\ad_\t C(\tau)$, which is the 
expected expression for $dt/d\tau$, and \refeq{phi_cst} boils down to 
$\phi(\xi,\tau)=C(\tau)$.
\end{proof}

 
\begin{prop}\label{prop_dA}
The following formula holds:
$$\frac{\partial A}{\partial\tau} = [\phi(\tau), A] + \psi_A$$
\end{prop}
\begin{proof}
We'll get this by asserting the flatness of $\nabla$. The calculations are 
parallel to the proof of flatness of $\nabla_\PC$ of \ref{sec_univ}.

The 1-dimension assumption gives us $d(\phi(\tau)d\,\tau) = 0$. Let's compute 
the other terms.
\begin{eqnarray*}
d\nu &=& -\left(\frac{\partial \t}{\partial \tau}\frac{\partial }{\partial \t} 
         +\frac{\partial A}{\partial \tau}\frac{\partial }{\partial A}\right)
         \nu \wedge d\,\tau\\
 &=& \left(\left(-\frac{A}{2\pi i} +[\phi, \t]\right) 
           \frac{\partial}{\partial\t}+\frac{\partial A}{\partial\tau}
           \frac{\partial}{\partial A}
      \right)
      \ad_\t F(\xi , \ad_\t ;\tau ) A  d\xi\wedge d\tau\\
 &=& \frac{(X+Y) F(\xi , X+Y ;\tau )- Y F(\xi ,  Y;\tau )}{X} 
     \left(\left(-\frac{A}{2\pi i} +[\phi, \t]\right),A\right)
     d\xi\wedge d\tau \\
   && + \ad_\t F(\xi , \ad_\t ;\tau )\frac{\partial A}{\partial \tau}
        d\xi\wedge d\tau
\end{eqnarray*}

The second summand, $\nu\wedge\nu$, equals, up to the factor 
$2^{-1}(2\pi i)^{-1}d\xi\wedge d\tau$,
$$\left(XYF(\xi,X+Y ;\tau)  \left(\wp(X) -\wp(Y) \right) +
\frac{XF(\xi , X  ;\tau)}{Y}-\frac{YF(\xi , Y  ;\tau )}{X}\right)\db{A,A}$$

The last term, 
$\nu\wedge (\phi(\tau)d\,\tau)+ (\phi(\tau)d\,\tau)\wedge \nu$, is equal to
$$[\ad_\t F(\xi , \ad_\t ;\tau )A,\phi(\tau)]d\xi
\wedge d\tau=XF(\xi , X;\tau )\db{A,\phi(\tau)}
d\xi \wedge d\tau.
$$

Finally, up to factor $\frac12\frac{d\,\xi\wedge d\,\tau}{2\pi i}$, the 
curvature is equal to:
\begin{eqnarray*}
&&  \ad_\t F(\xi , \ad_\t ;\tau )\left(
     \frac{XY}{X+Y}
     \left(\left(\wp(X)-\frac1{X^2}\right)-\left(\wp(Y)-\frac1{Y^2}\right)
     \right)\db{A,A}\right. \\
&& \qquad
   \left. + 4\pi i\left(\frac{\partial A}{\partial \tau}-\db{\phi(\tau),A}
                   \right)
     \right)\\
&=& \ad_\t F(\xi , \ad_\t ;\tau )
     \left(4\pi i\psi_A  + 
           4\pi i \left(\frac{\partial A}{\partial \tau}-\db{\phi(\tau),A}
                  \right)
     \right).
\end{eqnarray*}

By invertibility of $\ad_\t F(\xi , \ad_\t ;\tau )$, the vanishing of this 
curvature thus yields the wished formula.
\end{proof}

\sss{Conclusion for the local case}\label{local_ess_surj2}
Let's consider on $\VC$ the connection: 
$$\nabla_\VC\ \ass\ d-\phi(\tau) d\tau$$
The operators $\t(\tau)$ and $A(\tau)$ give the $\RC$-algebra structure on 
$\VC$ and propositions \ref{prop_phi_cst} and  \ref{prop_dA} 
state its horizontality. Applying $\FC_U$, we obtain the 
original $\VG$, thanks to formula \refeq{nabla_rel_ess_surj}.

In the higher dimensional case, there are more factors to take into account in 
the computations above. It turns out that they express nothing more 
than the flatness of $\nabla_\VC$. We'll spare the reader those details.

As for fully faithfullness, the proof in \ref{fw_full_faith} carries over 
transparently.

\sss{General bases}\label{sss_rel_glue}
According to the previous paragraph, there is an open covering
$(U_i)_{i\in I}$ of $S$ such that each $\FC_{U_i}$ is an
equivalence. Then one can construct an essential inverse $u$ to $\FC$ in the 
following way: 

Let us choose an essential inverse $u_i$ of each $\FC_{U_i}$. The
functors $u_i|U_i\cap U_j$ and $u_j|U_i\cap U_j$ are both inverses of
$\FC_{U_i\cap U_j}$ and therefore canonically isomorphic through
$\phi_{ij}$. The canonicity implies in particular that the descent
condition $\phi_{ik}=\phi_{ij}\phi_{jk}$ holds for any triple
$\{i,j,k\}$  
For any object $\VG$ of $\nilpconn(X/S;e)$, the
various $u_i(\VG)$ therefore define an object $u(\VG)$ of
$\mods{\RC}$, and it goes in the same way for morphims, which can be
defined locally. It's easy to see that $u$ and $\FC$ are inverse of
each others, again because arrows can be defined locally in both categories.
\hfill\qed

\ssect{The Fundamental Hopf algebras for geometrical fiber functors.}
\label{ss_geom}
We now consider a fixed family $X\to S$ and a section $\sigma$ of
$X\setminus e(S)\to S$. The pull-back $\sigma^*$ is a 
fiber functor of $\nilpconn(X/S; e)$ over $\conn(S)$, fixed by the 
canonical isomorphism $\s^*p^*\simeq\id$. In this subsection, we want to
identify the corresponding fundamental Hopf algebra in
the target category $\conn(S)$ at $\sigma^*$.

Since the base $S$ is fixed, we will silently drop it from notations.

\sss{}It is also worthwhile to consider 
the following generalization of such pull-backs, known as 
tangential base points. Let $z$ be a $1$-germ of transversal coordinate at 
the section $e$. Then it defines a fiber functor as follows: 
let  $\nabla$ be a flat 
connection on $X$ with simple pole along $e$; denote by $K$ its
residue at $e$. Then the local connection $\nabla-Kdz/z$ has no
singularities and it is simple to check that
$e_z^*(\nabla):=e^*(\nabla-Kdz/z)$ is a flat connection on $S$.
Note that this functor depends on the choice of $1$-germ of $z$. 
A more canonical version would involve the relative normal bundle 
of $e$ in $X$.

\sss{Twistings}
As was earlier mentioned, a continuation of the fiberwise functors 
${\cal F_s}\colon\mods{R}\to\nilpconn(X_s;e), s\in S$
  to the relative functor ${\cal F}_S\colon
\mods{{\cal R}} \to \nilpconn(X/S;e)$ is not unique. Indeed, let
$\mu$ be a $\RC$-valued differential $1$-form on $S$
that satisfies the Maurer-Cartan equation:
\begin{equation}
  \def\theequation{MC}
  \nabla_{\!\RC}\mu + \mu\wedge\mu = 0 
\end{equation}
Then we can define a new connection 
$\nabla_{\cal R}^\mu\ \ass\  \nabla_{\cal R}+\ad_{\mu}$ 
on the $S$-bundle ${\cal R}$. It is easy to check 
that (MC) yields the flatness of $\nabla_{\cal R}^\mu$. 
We'll denote by ${\cal R}^\mu$ the bundle $\RC$ equipped with this connection.
Further, if $\mu$ takes values in primitive elements of $\RC$, 
then $\RC^\mu$ is again a Hopf algebra in $\conn(S)$, 

Analogously, one can also define a new flat connection on ${\frak P}$ as 
$\nabla_{\PG}^\mu\ \ass\ \nabla_{\PG} - p^*\mu$, where $\mu$ acts from the 
right. This bundle with connection
will be denoted by ${\PG}^\mu$. It is a 
$p^*\RG^\mu$ right module and induces a functor
$${\cal F}^\mu_S\colon\mods{{\cal R}^\mu} \to \nilpconn(X/S;e)$$ 
by the same rule \refeq{def_F} that ${\frak P}_{X/S}$ induces ${\cal F}_S$.

\sss{Twisted equivalence}
 Note that the categories $\mods{{\cal R}}$ and $\mods{{\cal R}^\mu}$
 are equivalent under the functor $F_{0,\mu}$ sending 
 $({\cal V},\nabla_{\cal V})$ to $({\cal V},\nabla_{\cal V}+\mu)$. 
 In Morita form, this functor is given by the
 bimodule ${\cal P}_{0,\mu}$ of $\bimods{{\cal R}^\mu}{\cal R}$, 
 which is defined as the bundle ${\cal R}$ 
 with connection $\nabla_{\cal R}+l_\mu$ (left multiplication).

 Evidently, ${\cal F}$ is isomorphic to ${\cal F}^\mu\circ F_{0,\mu}$, 
 so ${\cal F}^\mu_S$ is an equivalence of categories. This defines a 
 new fiber functor $\omega^\mu_\dr$, again fixed by $\s^*p^*\simeq\id$. 
 We can think of
 ${\cal P}_{0,\mu}$ as the torsor of paths from $\omega_\dr$ to 
 $\omega^\mu_\dr$. More generally, the torsor of paths from 
 $\omega^\mu_\dr$ to $\omega^{\mu'}_\dr$
 would be the bimodule
 ${\cal P}_{\mu,\mu'}$ in $\bimods{{\cal R}^{\mu'}}{\cal R}^\mu$, whose 
 connection is $\nabla_{\cal  R}+l_{\mu'}-r_\mu$ 
 (left and right multiplications).

\sss{Equivalence at $\sigma$}
 We now turn to the composition $\sigma^*\circ {\cal F}^\mu_S$.
 By the very definition, we get that it maps $({\cal V},
 \nabla_{\cal V})$ to $({\cal V}, \nabla_{\cal
 V}+\nu_\sigma-\mu)$, where $\nu_\sigma$ is defined as 
 $\nabla_\VC-p^*\nabla_\PG$. More precisely, this is true up to 
 $\s^*p^*\simeq\id$.  

 In the local $S=\HM$ situation, this $\nu_\sigma$ 
 is nothing but $\s^*\nu$, with $\nu$ as in \refeq{def_P_rel}. 
 This motivates the notation and shows that 
 $\nu_\sigma$ is actually independent of $\nabla_\VC$.

 For $\mu=\nu_\sigma$, we obtain that the fixing isomorphism $\s^*p^*\simeq\id$
 induces a canonical isomorphism from this composition to 
 the forgetful functor from $\mods{{\cal R}}^\mu_S$ to $\conn(S)$. Together 
 with the fact that ${\cal F}^\mu_S$ 
 is an equivalence of categories, this allows us to conclude that 
 ${\cal R}_S^{\nu_\sigma}$ is the fundamental Hopf algebra at the fiber functor
 $\sigma^*$, fixed by $\s^*p^*\simeq Id$. 
 In the same way, for two sections $\sigma_1$ and $\sigma_2$,
 the bimodule ${\cal P}_{\mu_1,\mu_2}$ with $\mu_i=\nu_{\sigma_i}$ is
 the torsor of paths from $\sigma_1$ to $\sigma_2$.

\sss{A differential equation for parallel transports}
Since the Riemann-Hilbert correspondence commutes with geometric functors
like $p^*$, $\s^*$ and respects natural isomorphisms as $\s^*p^*\simeq\id$
we conclude from the previous paragraph that the parallel transport $\Phi$ 
with respect to the connection
${\frak P}_{X/S}$ along a $S$-family of paths from $\sigma_1$ to
$\sigma_2$ satisfies the equation
$$\nabla \Phi + \nu_{\sigma_1}\Phi - \Phi\nu_{\sigma_2}=0.$$

Let us stress that for $\HM=S$, it takes the following very explicit form:
$$\nabla \Phi + \sigma_1^*(\nu)\Phi-\Phi\sigma_2^*(\nu)=0.$$

This would remain true in the case of tangential base points.

\section{Algebraic $\QM$-structure}\label{sec_alg}
In this section, we exhibit the natural algebraic $\QM$-structures of the 
previously constructed fundamental groups and torsors. By ``natural'', we mean 
here with respect to the fact \cite[10.41]{DelPi1} that the algebraic De Rham 
fundamental group(oid) 
(defined by means of algebraic connections over the base 
field) commutes with the extension of scalars to $\CM$. 
We should stress that the nilpotency 
assumptions we've made all along are crucial for this to be true \cite[10.35]{DelPi1}.

\ssect{Over a single curve}\label{ss_Qsingle}
We first treat the algebraic $\Bbb Q$-structures on $\cal P$ and $R$ 
for a single elliptic curve defined over $\Bbb Q$.

\sss{Reminders}
Let us remind some standard facts about algebraic elliptic curves, 
mostly to fix terminology and notations.

Let $(X,O)$ be an elliptic curve defined over some field $k$ of characteristic 
zero.
 Then it can be represented as a plane cubic $y^2=4x^3-g_2 x -g_3$,
 with $g_2, g_3\in k$ and the
 marked point $O$ lies at infinity. This form is unique up to
 dilatations $g_2\to\lambda^4g_2,\quad g_3\to\lambda^6g_3,\quad
 x\to\lambda^2x,\quad y\to\lambda^3y$, for $\lambda\in k$.

  If $k$ is a subfield of $\Bbb C$, e.g.,
$k=\Bbb Q$, one can represent the analytic curve $X(\Bbb C)$ as the
quotient of $\Bbb C$ by some lattice $\Lambda$. The
map  from $X(\Bbb C)$ to $\Bbb C/\Lambda$ is defined by integration
of the form
$\omega = dx/y$ from the marked point to the variable one
$p\to\int_O^p\omega$. The ambiguity in choosing of the path of integration
belongs to the topological first homology group $H_1(X(\Bbb C),\Bbb Z)$
which corresponds to the lattice $\Lambda$.

Choose some basis $(u, v)$ of $\Lambda$ in
such a way that $\Im (v/u)>0$, and put $\tau=v/u$. Then we have
$g_2=60u^{-4}e_4(\tau)$ and $g_3=140u^{-4}e_6(\tau)$. Moreover,
the map
$$
\xi\to (x(\xi)=u^{-2}(E_2(\xi,\tau)-e_2(\tau)),
\quad y(\xi)=-2u^{-3}E_3(\xi,\tau))$$
produces an analytic isomorphism between $X_{\tau}=\Bbb C/(\Bbb Z\tau +\Bbb Z)$
and  $X(\Bbb C)$ and we have $\omega=dx/y=ud\xi$. So, 
$u$ can be retrieved as the
period (elliptic integral) of the algebraic form $\omega$ against a
topological chain. 
If we choose another equation of $X$, twisted, say, by $\lambda$,
then $\omega$ becomes $\lambda^{-1}\omega$ and $(u,v)\to (\lambda^{-1}u,
\lambda^{-1}v)$.
\sss{Some analytic preparations} Note that over the complement
$U_\tau$ of the marked point in $X_\tau$, the bundle $\cal P$ with multiplier $\exp(-2\pi i\t)$ of section \ref{sec_fw} is trivialized by the 
left multiplication by $g(\xi)=\exp(-E_1(\xi,\tau)\t)$. 

Indeed, sections of $\cal P$ are the quasiperiodic functions $f$ on $\Bbb C$:
$f(\xi+1)=f(\xi);\quad f(\xi+\tau)=\exp(-2\pi i\t)f(\xi)$. For each such, 
$g(\xi)f$ is elliptic, since we have $E_1(\xi+1,\tau)=E_1(\xi,\tau)$ and 
$E_1(\xi+\tau,\tau)= E_1(\xi,\tau)-2\pi i$.

\sss{}The connection $\nabla=d-\ad_\t F(\xi,\ad_\t)Ad\xi$
transforms under the gauge transformation $s\to g(\xi)s$ into the operator
$$\nabla_\alg\ \ass\ 
g(\xi)\nabla g^{-1}(\xi) = d-dg(\xi)g^{-1}(\xi)-
g(\xi)\ad_\t F(\xi,\ad_\t)Ag^{-1}d\xi,$$
which we expand in terms of Eisenstein series using \refeq{F_Eis}:
\begin{gather*}
= d - E_2(\xi,\tau)\t d\xi
 -\exp(-E_1(\xi,\tau)\ad_\t )
 \exp\left(-\sum _{k=1}^{\infty}\frac{(-\ad_\t)^k}k
(E_k(\xi,\tau  )-e_k(\tau))\right)A\,d\xi \\
=d-\left(
(E_2(\xi,\tau)-e_2(\tau))\t    +
 \exp\left(-\sum _{k=2}^{\infty}\frac{(-\ad_\t  )^k}k
(E_k(\xi,\tau)-e_k(\tau))\right)(A+e_2(\tau)\t)\right)d\xi.
\end{gather*}

\sss{The $\QM$-structure}
Let us introduce the following new generators of the ring $R$:
$$\talg=u\t,\qquad\salg=u^{-1}(A+e_2(\tau)\t)$$
In terms of these, $\nabla_\alg$ takes the form $d-K(\xi,\tau)ud\xi$, 
where $K(\xi,\tau)$ is:
$$ u^{-1}(E_2(\xi,\tau)-e_2(\tau))\talg    +
 \exp\left(-\sum _{k=2}^{\infty}\frac{(-\ad_{\talg } )^k}k
u^{-k}(E_k(\xi,\tau  )-e_k(\tau))\right){\salg}
$$

As the $e_{2k}(\tau)$, $k\ge 4$ are polynomials in $e_4(\tau)$ and
$e_6(\tau)$ with rational coefficients, the
$u^{-2k}e_{2k}(\tau)$  are polynomials in $60u^{-4}e_4(\tau)=g_2$ and
$140u^{-6}e_6(\tau)=g_3$ with rational coefficients.

As the $E_k(\xi,\tau)$ are polynomials in $E_2(\xi,\tau)-e_2(\tau)$,
$E_3(\xi,\tau)$, $e_4(\tau)$ and $e_6(\tau)$ with rational coefficients,
the $u^{-k}E_k(\xi,\tau)$ are polynomials in 
$u^{-2}(E_2(\xi,\tau)-e_2(\tau)=x$,
$-2u^{-3}E_3(\xi,\tau)=y$, $g_2$ and $g_3$ with rational coefficients.

So, we see that for generators $\talg$ and ${\salg}$,
the transformed connection $\nabla_{\alg}$ is algebraic over $\Bbb
Q$; hence the ring $R_{alg}=\sernc{\QM}{\talg, \salg}$
is the natural algebraic $\Bbb Q$-structure on $R$.  Note that this
rational structure does not depend in the choice of the equation
of the elliptic curve, as a change of the equation multiplies $u$
by some {\it rational} number $\lambda^{-1}$.
\joenote{Remark about apparent irregularity?}
\ssect{Relative case}
We start by a short account in our notations of the realisation of the 
analytic stacks $\Bbb H/SL_2(\Bbb Z)$ and $\Bbb X/SL_2(\Bbb Z)$
as stacks of $\Bbb C$-points of some algebraic stacks. All the results 
mentionned there are very classical and belong to the mathematics of 19{th} 
century. We then apply the same gauge transformation as for single elliptic 
curves.
\sss{}
Consider the product $\Bbb H\times \Bbb C^*$ and 
denote the coordinate on $\Bbb C^*$
by $u$. Define the action of $SL_2(\Bbb Z)$ by the following formula:
$$\left(\begin{array}{cc}
a&b\\
c&d
\end{array}\right):(\tau,u)\to\left(\frac{a\tau+b}{c\tau+b},
(c\tau+d)u\right).$$ 
The map
$(\tau,u)\to (g_2=60u^{-4}e_4(\tau),\quad g_3=140u^{-4}e_6(\tau))$
is $SL_2(\Bbb Z)$-invariant and provides an isomorphism of the
quotient and  the complement ${\cal B}_{\Bbb C}$ of the curve
$\Delta:=g_2^3-27g_3^2=0$ in $\Bbb C^2$. As $ \Bbb H/SL_2(\Bbb Z)$
is the quotient of $ (\Bbb H\times\Bbb C^*)/SL_2(\Bbb Z)$ by $\Bbb
C^*$ acting on the second factor, this quotient stack  is equal to
$\Bbb C$-points of the quotient stack of ${\cal B}=\Bbb
A^2\setminus \{\Delta =0\}$ by the following action of $\Bbb G_m$:
$(g_2, g_3)\to(\lambda^{-4}g_2, \lambda^{-6}g_3)$. Evidently, this
algebraic stack is defined over $\Bbb Q$.

In the same way, we can consider the quotient $(\Bbb X\times\Bbb C^*)/SL_2(\Bbb Z)$ with
the same action of the group $SL_2(\Bbb Z)$ on the second factor.\joenote{Isn't it $\Gm$?}
Furthermore, the map $(\tau,\xi,u)\to  (g_2=60u^{-4}e_4(\tau),\quad g_3=140u^{-4}e_6(\tau)
x =u^{-2}(E_2(\xi,\tau)-e_2(\tau)),y=-2u^{-3}E_3(\xi,\tau),z=1)$ is an
isomorphism of this quotient with the complex projective cubic ${\cal Q}_{\Bbb C}
=\{y^2z=4x^3-g_2xz^2-g_3z^3\}$
over ${\cal B}_\Bbb C$. 

We get that $\Bbb X/SL_2(\Bbb Z)$ is the stack of $\Bbb
C$-points of the quotient of the cubic ${\cal Q}=\{y^2z=4x^3-g_2xz^2-g_3z^3\}$
over $\cal B$ by the action of $\Bbb G_m$ defined as follows:
$$(g_2,g_3,x,y)\To{\lambda}(\lambda^{-4}g_2, \lambda^{-6}g_3 ,\lambda^{-2}x,
\lambda^{-3}y)$$ 
The variety $\cal Q$ and its $\Gm$-action on it are defined over
$\Bbb Q$. The complement of the neutral section corresponds to the
affine curve ${\cal Q}^{\aff}=\{y^2=4x^3-g_2x-g_3\}$.

\sss{Analytic preparation} As in the case of individual curve we apply the 
gauge transformation of the left multiplication by 
$g(\xi,\tau)=\exp(-E_1(\xi, \tau)\t)$. This transforms the connection
$\nabla$ on $\PG$ in
\begin{eqnarray*}
\nabla_\alg  &=& d-dg(\xi,\tau)g^{-1}(\xi,\tau)
                  +g(\xi,\tau)\nu(\xi,\tau)g^{-1}(\xi,\tau) \\
              &&  -g(\xi,\tau)\psi_\t\partial_\t g^{-1}(\xi,\tau)
                  -g(\xi,\tau)\psi_A\partial_Ag^{-1}(\xi,\tau)
\end{eqnarray*}
\begin{gather*}
=d-\left(E_2 d\xi -\frac1{2\pi i}(E_3-E_1E_2)d\tau\right)\t
-\left(\exp\left(-\sum _{k=2}^{\infty}\frac{(-\ad_\t )^k}k
(E_k -e_k )\right) d\xi\right.\\
+\left(\left(-\frac{1}{\ad_\t}
+\sum _{k=1}^{\infty}
(-\ad_\t )^{k-1}
(E_k -e_k )
\right)
\exp\left(-\sum _{k=2}^{\infty}
\frac{(-\ad_\t )^k}k
(E_k -e_k )\right)\right.\\
\left.\left.+
\frac{\exp(-E_1\ad_\t)}{\ad_\t}\right)\frac{d\tau}{2\pi i}\right) A
+\left(\psi_\t\partial_\t -\frac1{2\pi
i}\frac{\exp(-E_1\ad_\t)-1}{-\ad_{\t}}
Ad\tau\right)+\psi_A\partial_A\\
=d-\left(E_2 \t
+\exp\left(-\sum _{k=2}^{\infty}\frac{(-\ad_\t )^k}k
(E_k -e_k )\right)A \right)\left(d\,\xi +\frac1{2\pi i}E_1d\tau\right)
+ E_3 \frac{\t d\,\tau}{2\pi i} \\
-\left(\left(-\frac{1}{\ad_\t}
+\sum _{k=2}^{\infty}
(-\ad_\t )^{k-1}
(E_k -e_k )
\right)
\exp\left(-\sum _{k=2}^{\infty}
\frac{(-\ad_\bfa )^k}k
(E_k -e_k )\right)+
\frac{1}{\ad_\t}\right)\frac{Ad\,\tau}{2\pi i}\\
+\psi_\t\partial_\t  +\psi_A\partial_A\\=
d+\nu_\alg+\psi_\t\partial_\t  +\psi_A\partial_A.
\end{gather*}

In this expansion, we used the equality
$2\pi i \partial E_1/\partial\tau= E_3-E_1E_2$.

If we replace generators $\t, A$ by $\talg=u\t, {\salg}=u^{-1}(B+e_2\t)$, 
the coefficients of the differential form $\nu_{\alg}$
become polynomials in $u^{-2}(E_2-e_2)=x$, $-2u^{-3}E_3=y$, $60u^{-4}e_4=g_2$ 
and $140u^{-6}e_6=g_3$. What remains to be done for $\nu_\alg$ is summarized as
the following proposition.

\begin{prop}The differential forms 
$$\frac{u^2d\,\tau}{2\pi i},\quad \left(e_2\frac{d\,\tau}{2\pi
i}+\frac{d\,u}u\right), \quad \mbox{and}\quad
u^{-1}\left(d\xi +\frac1{2\pi i}E_1d\tau\right)$$
can be expressed as rational polynomials in terms of $x, y, z, g_2$ and $g_3$.
\end{prop}

\begin{proof}Note that
$d\Delta=-12 \Delta (e_2(2\pi i)^{-1}d\tau+u^{-1}du)$, so
$e_2(2\pi i)^{-1}d\tau+u^{-1}du$
is the algebraic differential form $\kappa\ \ass\ 12^{-1}\Delta^{-1}d\Delta$.
We can treat the operator $d + j \kappa$ as a $\Bbb G_m$-equivariant 
connection on the trivial bundle on $\cal B$, equipped with the action 
$\lambda\to\lambda^{-j}$ of $\Bbb G_m$. This explains
the coefficient of $\kappa$ in formulas below.

From $2\pi ide_4=(14e_6-4e_2e_4)d\tau$ follows
$dg_2=6g_3(2\pi i)^{-1}u^2d\tau-4\kappa g_2$, hence
$(2\pi i)^{-1}u^2d\tau$ 
is algebraic.

As $2\pi i\partial e_2/\partial\tau=(5e_4-e_2^2) $ and
$2\pi i\partial E_2/\partial\tau=3((E_2-e_2)^2-5e_4)-E_2^2-2E_1E_3$, we have
$$  dx=y\left(u(d\,\xi+\frac1{2\pi i}E_1d\,\tau)\right)
+\left(2x^2-\frac13g_2\right)
\frac{u^2d\tau}{2\pi i}-2x\kappa ,
$$
\end{proof}
\sss{}
As we perform the change of generators $(\t, A)\to (\talg, {\salg})$
the differential term $d+\psi_\t\partial_\t + \psi_A\partial_A$
transforms.
The generator $\salg=u^{-1}(A+e_2\t)$ depends directly in $\tau$, and both  
$\talg=u\t$ and $\salg$ depend in $u$, so the differentiation $d$
transforms into 
$d+\talg\partial_{\talg}u^{-1}du+\talg\partial_{\salg}u^{-2}de_2
-\salg\partial_{\salg}u^{-1}du$; furthermore
$\psi_\t\partial_\t=-\frac{Ad\,\tau}{2\pi i}\partial_{\t}$ transforms into
$$\left(u^{-1}e_2\talg - u{\salg}\right) \frac{d\,\tau}{2\pi i}
  \left(u\partial_{\talg} 
+ u^{-1}e_2\partial_{\salg}\right).$$
As $de_2=(2\pi i)^{-1}(5e_4-e_2^2)d\tau$, the term 
$d + \psi_\bfa\partial_\bfa$ becomes:

$$d+  \left(\kappa\talg +\frac{u^2d\,\tau}{2\pi i}\salg\right)\partial_{\talg}
+\left( \frac1{12}g_2\frac{u^2d\,\tau}{2\pi i}\talg
-\kappa\salg\right)\partial_{\salg};$$
which is $\Bbb Q$-algebraic.

Finally, the summand $\psi_A\partial_A$ transforms into
$$\psi_{\salg}\partial_{\salg}=-\frac12
u\frac{u^{-2}XY}{u^{-1}X+u^{-1}Y}
$$
$$\times\left(\left(\wp(u^{-1}X)-\frac1{(u^{-1}X)^2}\right)
-\left(\wp(u^{-1}Y)-\frac1{(u^{-1}Y)^2}\right)
\right)\odb {\salg},{\salg}\cdb_{\talg} \frac{d\,\tau}{2\pi i}\partial_{\salg},$$
and the coefficient of each term
$[\ad_{\talg}^i{\salg},\ad_{\talg}^j{\salg}]\partial_{\salg}$
 in this expression can be calculated in terms of  the Taylor
 expansion at zero of $\wp(X)-X^2$. This coefficient
 is equal to $(-1)^{j-1}(i+j)(2\pi i)^{-1}u^{i+j+1}e_{i+j+1}u^2d\tau$, 
 and is therefore a polynomials in $g_2$ and $g_3$ 
 multiplied by the algebraic form $(2\pi i)^{-1}u^2d\tau$.

So we constructed $\Bbb Q$-algebraic connections 
$({\cal R}_{\alg}, \nabla_{{\cal R}_\alg})$ on ${\cal B}$ and  
$({\frak P}_{\alg}, \nabla_{{\frak P}\alg})$ on
 ${\cal Q}^{\aff}$ which are equal to the analytic ones after tensoring by 
$\Bbb C$. These algebraic connections are moreover $\Bbb G_m$-equivariant.  

\sss{Pull-back to a family}
 Let $X\to S$ be a smooth algebraic family of elliptic curves defined over 
 $\Bbb Q$. Then
we have a well defined map from the family $X(\Bbb C)\to B(\Bbb C)$
of $\Bbb C$-points to ${\Bbb X}/SL_2(\Bbb Z)$, and this map determines
the analytic connections ${\cal R}_S$ and ${\frak P}_{X/S}$, on which
we shall now define an algebraic $\Bbb Q$-structure.

Denote the generic point of $S$ by $\eta$. 
The elliptic curve $X_\eta$. Since $\eta$ is the spectrum of the field 
$\QM(\eta)$, the elliptic curve $X_\eta$ over $\eta$ is isomorphic to a 
plane cubic $y^2=4x^3-g_2 x -g_3$, whose equation is unique up to
 dilatations $g_2\to\lambda^4g_2,\quad g_3\to\lambda^6g_3,\quad
 x\to\lambda^2x,\quad y\to\lambda^3y$, with $\lambda\in \QM(\eta)$. 

 Fix some choice of $g_2$ and $g_3$, let $\tilde{S}$ be the open subscheme 
 of $S$ determined by inequalities $g_2\neq \infty$ , $g_3\neq \infty$ and $\Delta\neq 0$, and consider the induced family $\tilde{X}\to \tilde{S}$. 
We have a map
 $g=(g_2,g_3)$ from $\tilde{S}$ to $\cal B$ and $\tilde{X}$ is induced by $g$.
 
 The connections ${\cal R}_{\tilde{S}}$ and ${\frak P}_{\tilde{X}/\tilde{S}}$
 are equal to the pull-back of the corresponding connections on ${\cal
 Q}\to{\cal B}$, so the pull-back of ${\cal R}_{\alg}$ and ${\frak P}_{\alg}$
 provide a $\Bbb Q$-algebraic structure on them. Note that  
  $g^*{\cal R}_{alg}$ and $g^*{\frak P}_{alg}$ can be continued  as regular\joenote{regular singular}
 algebraic connections to $S$ and $X$ respectivly, as after tensoring by 
 $\Bbb C$ they become regular analytic connections.
 
 These algebraic $\Bbb Q$-structures don't depend in the choice of $g_2$ and
 $g_3$ as they\joenote{$g_2$ and $g_3$?} can be changed by dilated one and  
 ${\cal R}_{ {\cal B}}$ and ${\frak P}_{{\cal
 Q}/{\cal B}}$ are $\Bbb G_m$-equivariant.

\sss{Twistings}Note that if $X/S$ be an algebraic family over $\Bbb Q$ and 
$\mu$
is $\Bbb Q$-algebraic ${\cal R}_{alg}$-valued $1$-form, then
${\cal R}^\mu$ and ${\frak P}^\mu$ can be equipped by an algebraic
$\QM$-structure. Hence, for a $\sigma$, defined over $\QM$,  the fundamental
Hopf algebra at $\sigma$ has also a natural algebraic $\Bbb Q$-structure.

\bibliographystyle{joeplain}
\bibliography{defs,geometrie,traites,motifs,Dmodules,polychoses,divers,moi,ICM}

\def\noop#1{}\def\noop#1{}
\providecommand{\bysame}{\leavevmode ---\ }
\providecommand{\og}{``}
\providecommand{\fg}{''}
\providecommand{\smfandname}{et}
\providecommand{\smfedsname}{\'eds.}
\providecommand{\smfedname}{\'ed.}
\providecommand{\smfmastersthesisname}{M\'emoire}
\providecommand{\smfphdthesisname}{Th\`ese}
\begin{thebibliography}{10}

\bibitem{Andre}
{\scshape Y.~Andr{\'e}} -- {\og Diff\'erentielles non commutatives et th\'eorie
  de {G}alois diff\'erentielle ou aux diff\'erences\fg}, \emph{Ann. Sci.
  \'Ecole Norm. Sup. (4)} \textbf{34} (2001), no.~5,
p.~685--739.

\bibitem{Atiyah_fibres}
{\scshape M.~F. Atiyah} -- {\og Vector bundles over an elliptic curve\fg},
  \emph{Proc. London Math. Soc. (3)} \textbf{7} (1957),
p.~414--452.

\bibitem{BeilDel}
{\scshape A.~Be{\u\i}linson {\normalfont \smfandname} P.~Deligne} -- {\og
  Interpr\'etation motivique de la conjecture de {Z}agier reliant
  polylogarithmes et r\'egulateurs\fg}, in \emph{Motives} \cite{MOT},
  p.~97--121.

\bibitem{BeilLev}
{\scshape A.~Be{\u\i}linson {\normalfont \smfandname} A.~Levin} -- {\og The
  elliptic polylogarithm\fg}, in \emph{Motives} \cite{MOT}, p.~123--190.

\bibitem{BoAl}
{\scshape A.~Borel, P.-P. Grivel, B.~Kaup, A.~Haefliger, B.~Malgrange
  {\normalfont \smfandname} F.~Ehlers} -- \emph{Algebraic ${D}$-modules},
  Persp. in Math., no.~2, Academic Press, Boston, 1987.

\bibitem{Del_singreg}
{\scshape P.~Deligne} -- \emph{\'{E}quations diff\'erentielles \`a points
  singuliers r\'eguliers}, Lect. Notes in Math., no. 163, Springer-Verlag,
  Berlin, 1970.

\bibitem{DelPi1}
\bysame , {\og Le groupe fondamental de la droite projective moins trois
  points\fg}, Galois groups over {$\QM$} (Berkeley, CA, 1987), Springer, New
  York, 1989, p.~79--297.

\bibitem{DelFest}
\bysame , {\og Cat{\'e}gories tannakiennes\fg}, The Grothendieck Festschrift,
  Vol.\ II, Birkh{\"a}user, Boston, MA, 1990, p.~111--195.

\bibitem{DelGon}
{\scshape P.~Deligne {\normalfont \smfandname} A.~B. Goncharov} -- {\og Groupes
  fondamentaux motiviques de tate mixte\fg}, Publication {\'e}lectronique,
f{\'e}vrier 2003, \arxiv{math.NT/0302267}.

\bibitem{ENR}
{\scshape M.~Espie, J.-C. Novelli {\normalfont \smfandname} G.~Racinet} -- {\og
  Formal computations about multiple zeta values\fg}, \emph{From Combinatorics
  to dynamical systems, Journ{\'e}es de calcul formel, Strasbourg 2002}
  (F.~Fauvet {\normalfont \smfandname} C.~Mitschi, \smfedsname), IRMA Lectures
  in Mathematics and Theoretical Physics, vol.~3, de Gruyter, 2003.

\bibitem{Gonch98}
{\scshape A.~B. Goncharov} -- {\og Multiple polylogarithms, cyclotomy and
  modular complexes\fg}, \emph{Mathematical Research Letters} \textbf{5}
  (1998),
p.~497--516.

\bibitem{Gonch2001}
\bysame , {\og {Multiple polylogarithms and mixed Tate motives}\fg},
  Publication {\'e}lectronique,
mars 2001, \arxiv{math.AG/0103059}.

\bibitem{IhICM}
{\scshape Y.~Ihara} -- {\og Braids, {G}alois groups, and some arithmetic
  functions\fg}, \emph{1990 Int. Cong. of Mathematicians}, Math. Soc. of
  Japan., Tokyo, 1991, p.~99--120.

\bibitem{MOT}
{\scshape U.~Jannsen, S.~Kleiman {\normalfont \smfandname} J.-P. Serre}
  (\smfedsname) -- \emph{Motives, proceedings of the {AMS}-{IMS}-{SIAM} joint
  summer research conference}, Univ. of Washington, Seattle 1991, American
  Math. Soc., 1994.

\bibitem{WildesMixedFal}
{\scshape {J{\"o}rg Wildeshaus}} -- {\og {Mixed structures on fundamental
  groups}\fg},
{Preprint, November 19, 1994, K-theory Preprint Archives,
  http://www.math.uiuc.edu/K-theory/0042/}.

\bibitem{KroneckerF}
{\scshape Kronecker} -- {\og Zur theorie der elliptischen {F}unktionen\fg},
  Mathematische Werke, vol.~IV, 1881, p.~313--318.

\bibitem{LevinAnal}
{\scshape A.~Levin} -- {\og Elliptic polylogarithmes: an analytic theory\fg},
  \emph{Compos. Math} \textbf{106} (1997), no.~3,
p.~267--282.

\bibitem{joeIHES}
{\scshape G.~Racinet} -- {\og Doubles m{\'e}langes des polylogarithmes
  multiples aux racines de l'unit{\'e}\fg}, \emph{Pub. Math. IH{\'E}S}
  \textbf{95} (2002),
p.~185--231, \arxiv{math.QA/0202142}.

\bibitem{Reut}
{\scshape C.~Reutenauer} -- \emph{Free {L}ie algebras}, London Math. Soc.
  Monographs, New series, no.~7, Oxford, 1993.

\bibitem{spitzdiplom}
{\scshape M.~Spitzweck} -- {\og {R}ealisierungen motivischer
  {F}undamentalgruppen und ihre {A}nwendung auf die {B}eschreibung der
  {P}olylogarithmusgarben auf ${\mathbb{p}}^1\setminus\{0,1,\infty\}$\fg},
  {D}iplomarbeit, {U}niversit{\"a}t Bonn, http://www.uni-math.gwdg.de/spitz/,
1998.

\bibitem{WeilEll}
{\scshape A.~Weil} -- \emph{Elliptic functions according to {E}isenstein and
  {K}ronecker}, Springer-Verlag, Berlin, 1976, Ergebnisse der Mathematik und
  ihrer Grenzgebiete, Band 88.

\bibitem{Wildes}
{\scshape J.~Wildeshaus} -- \emph{Realizations of polylogarithms}, Lect. Notes
  in Math., no. 1650, Springer, 1997.

\bibitem{ZagModTheta}
{\scshape D.~Zagier} -- {\og Periods of modular forms and {J}acobi theta
  functions\fg}, \emph{Invent. Math.} \textbf{104} (1991), no.~3,
p.~449--465.

\end{thebibliography}

\end{document}